\documentclass[namedate,webpdf,imanum]{ima-authoring-template}

\usepackage[T1]{fontenc}
\usepackage[utf8]{inputenc}
\usepackage{amsmath}
\usepackage{amssymb}
\usepackage{url}
\usepackage{verbatim}
\usepackage{listings}
\usepackage{color}
\usepackage{array}
\usepackage{graphicx}
\usepackage{overpic}
\usepackage{epstopdf}
\usepackage{upquote}
\usepackage{epsfig}
\usepackage{relsize} 
\usepackage{mathtools}
\usepackage{courier} 
\usepackage{xcolor} 
\usepackage{tikz} 
\usetikzlibrary{shapes,arrows}
\usepackage{soul}
\usepackage{float}
\usepackage{xcolor}
\usepackage{bm} 

\definecolor{gray}{rgb}{0.5,0.5,0.5}
\definecolor{dkgreen}{rgb}{.068,.578,.068}
\definecolor{dkpurple}{rgb}{.320,.064,.680}
\usepackage{upquote}
\usepackage{listings}
\definecolor{mygreen}{rgb}{0.2,0.4,0.2}
\definecolor{mymauve}{rgb}{0.58,0,0.82}
\lstdefinestyle{customc}{
  language=Matlab,
  showstringspaces=false,
  basicstyle=\footnotesize\ttfamily,
  keywordstyle=\color{blue},
  commentstyle=\color{dkgreen},
  identifierstyle=\color{black},
  stringstyle=\color{dkpurple},
  deletekeywords={sech, zeros, plot, eye, sum, isnan, diff, spy, hold, max, expm, sign, sqrt, diag, ans, format, surf, axis, subplot, pi, sin, cos, eps, schur, sinh, imag, abs, cosh, asinh, realmax, nargin, exp, length, end, cumsum, feval, norm, title, num2str, fprintf, figure, log, cot, linspace, csc, end},
  emph=[1]{end)},emphstyle=[1]\color{black}, 
  morekeywords={parfor},
  numbers=left,                    
  numbersep=10pt,
  numberstyle=\tiny\color{gray},
  xleftmargin=0in,
  xrightmargin=0in,
}
\usepackage{graphicx}
\usepackage{xcolor}

\usepackage{tabularx}
\usepackage{calc}

\newlength{\DepthReference}
\newlength{\HeightReference}
\newlength{\Width}%

\newcommand{\mycolorbox}[2][red]%
{%
    \settowidth{\Width}{#2}%
    \setlength{\fboxsep}{0pt}%
    {%
        \raisebox{0.15em}%
        {%
                \parbox[b][.5em][b]{\Width}{\centering \colorbox{#1}{\parbox[c][.6em][t]{\Width}{#2}}}%
        }%
    }%
}

\newcommand\be{\begin{equation}}
\newcommand\ee{\end{equation}}

\usepackage[left]{lineno}

\theoremstyle{thmstyletwo}%
\theoremstyle{definition}
\newtheorem{example}{Example}%
\newtheorem{remark}{Remark}%

\numberwithin{equation}{section}

\usepackage[numbered,framed]{matlab-prettifier}

\long\def\@tablecaption#1#2{%
	\begingroup%
	\fontsize{10bp}{12pt}\itshape\selectfont%
	\emph{#1}\enskip\enskip{#2\par}%
	\endgroup\vspace{\belowcaptionskip}}

\begin{document} 
  
 \DOI{DOI HERE}
\copyrightyear{2024}
\vol{00}   
\pubyear{2023}
\access{Advance Access Publication Date: Day Month Year}
\appnotes{Paper}
\copyrightstatement{Published by Oxford University Press on behalf of the Institute of Mathematics and its Applications. All rights reserved.}   
\firstpage{1}
 

\title[Spectral collocation for functional and  delay differential equations]{A spectral collocation method for functional and  delay differential equations}
 
\author{Nicholas Hale*
\address{\orgdiv{Department of Mathematical Sciences}, \orgname{Stellenbosch University}, \orgaddress{\street{Stellenbosch}, \postcode{7600},  \country{South Africa}}}}

\authormark{Nicholas Hale}

\corresp[*]{Corresponding author: \href{email:nickhale@sun.ac.za}{nickhale@sun.ac.za}}

\received{Date}{0}{Year} 
\revised{Date}{0}{Year}
\accepted{Date}{0}{Year}


\abstract{
A framework for Chebyshev spectral collocation methods for the numerical solution of functional and delay differential equations (FDEs and DDEs) is described. The framework combines interpolation via the barycentric resampling matrix with a multidomain approach used to resolve isolated discontinuities propagated by non-smooth initial data.
Geometric convergence in the number of degrees of freedom is demonstrated for several examples of linear and nonlinear FDEs and DDEs with various delay types, including discrete, proportional, continuous, and state-dependent delay.
The framework is a natural extension of standard spectral collocation methods and can be readily incorporated into existing spectral discretisations,  such as in Chebfun/Chebop, allowing the automated and efficient solution of a wide class of nonlinear functional and delay differential equations. 
}



\maketitle


\section{Introduction}\label{sec:intro}%
Delay differential equations (DDEs) arise in a variety of application areas, notably population dynamics~\citep{Gopalsamy1992, kuang1993}, epidemiology~\citep{Brauer2001,rihan2021}, control theory~\citep{kyrychko2010,michiels2007}, and more~\citep{erneux2009}. Below we review some basic concepts and refer the reader to any of the many excellent introductory texts on DDEs and FDEs for details~\citep{driver2012,hale1971,Kolmanovskii2007,michiels2007}).

The most elementary DDEs are first-order initial value problems (IVPs) of the form  
\begin{equation}\label{eqn:DDE}   
\begin{array}{lllll}
y' (t) &=& f (t, y(t), y(\tau(t))), \quad &   0 \le t \le T,\ \quad \tau(t) \le t,\\
y(t) &=& \phi(t),  &t \le 0,
\end{array}
\end{equation}
with $\phi$, $f$, and $\tau$ given. More complex models may involve 
multiple delays, 
systems of equations, higher-order derivatives, DDE boundary value problems (DDE-BVPs), DDE eigenvalue problems (DDE-EVPs), or partial delay differential equations (PDDEs). The most common types of delay are: 
discrete delay, where $\tau(t) = t - p$, $p > 0$; proportional (or pantograph-type) 
delay, where $\tau(t) = qt$, $0 < q < 1$; continuous 
delay, where $y'(t) = f(t, y(t), \int^t_0y(\tau)f(t-\tau)\,d \tau)$; 
or some collection/combination of the above. Furthermore, one may consider neutral-type DDEs, whereby $y'(t) = f (t, y(t), y(\tau(t)), y'(\tau(t)))$, or state-dependent delays, where $\tau = \tau(t, y(\tau))$. DDEs are a special case of {functional} differential equations (FDEs), which relax the condition that $\tau(t)\le t$.  
For simplicity, we present most of the methodology in this work in the context of DDEs. 
However, the assumption $\tau(t)\le t$ is not required in the framework we propose, making it 
 applicable for both DDEs and FDEs.

Like ordinary differential equations (ODEs), most DDEs and FDEs arising in practical settings do not have closed-form solutions, so some kind of approximation is required. Standard numerical methods for the approximation of DDEs like~(\ref{eqn:DDE})  typically involve Runge--Kutta-like solvers equipped with a suitable interpolation to evaluate the delay terms (the {\em method of steps}), or rewriting the DDE as a system of ODEs with additional variables to represent the delayed terms. For an overview of such methods, we refer the reader to texts such as~\citep{bellen2013,Engelborghs2000,Shampine2009}, and references therein. A major difficulty arising in DDEs and FDEs not present in ODEs is the propagation of discontinuities arising from non-smooth initial data. Locating and resolving these discontinuities is necessary in order to compute an accurate and efficient approximation, even for low order methods~\citep{baker1995}.  

For smooth ODEs, spectral collocation methods (also known as pseudospectral methods)  typically exhibit geometric convergence as the number of degrees of freedom is increased, making them particularly efficient when a high degree of accuracy is required~\citep{boyd2001, fornberg1996, trefethen2000}.  Moreover, collocation methods extend naturally to boundary value problems and differential-eigenvalue problems. 
Given the power of spectral methods in solving ODEs, it is natural to apply them to DDEs and FDEs, and there are numerous examples of this in the literature. We mention a few key papers here, but more detailed overviews of the development of spectral methods in these areas are given by~\cite{meng2019} and~\cite{jafari2021}.
   
One of the first spectral methods for DDEs was introduced by~\cite{fox1971}, where a Chebyshev-tau approach was used to solve the linear pantograph equation. Later, a Legendre-tau approach for DDEs with discrete delay was presented by~\cite{ito1991}. More recently there have been developments in using spectral methods to solve continuous delay DDEs~\citep{brunner2005,brunner2007}, nonlinear DDEs~\citep{ali2009}, neutral and state-dependent DDEs~\citep{Sedaghat2012}, DDE-EVPs~\citep{liu1989}, and PDDEs~\citep{engelborghs2001, butcher2011}. In addition, spectral collocation methods have proven useful for computing characteristic roots of DDEs~\citep{breda2005, wu2012}. Two existing MATLAB software packages closely related to that which we present here are Chebpack~\citep{trif2011, chebpack} and DDE-BIFTOOL~\citep{BIFTOOL, BIFTOOLMANUAL}. Chebpack implements a Chebyshev-tau approach allowing the ready solution of a variety of linear neutral-type DDEs with discrete and pantograph delay, including higher-order DDEs and DDE-EVPs. DDE-BIFTOOL, on the other hand, not only solves various forms of DDEs, but is also capable of computing extensive bifurcation analyses for dynamical systems with delay. For PDDEs and DDE-EVPs, DDE-BIFTOOL uses polynomial collocation at Gauss--Legendre nodes~\citep{engelborghs2001} in a similar way to that which we describe in Section~\ref{sec:main} below.

The present work describes a unifying framework based on Chebyshev spectral collocation and the barycentric resampling matrix (see Section~\ref{subsec:bary}) that can solve at once a variety and combination of the types of DDE and FDE problems discussed above. The aim is to achieve spectral accuracy 
using a discretisation that can be easily incorporated in Chebfun~\citep{Chebfun} to utilise its convenient syntax and efficient algorithms for nonlinear and adaptive solves. Rather than attempt to rigorously prove convergence results for some limited class of problems, we demonstrate by numerous examples the efficacy and generality of the proposed framework. To that end, the outline of this paper is as follows. In Section~\ref{sec:prelim} we provide an overview of some necessary preliminaries, specifically the barycentric interpolation formula and how it leads to the  pseudospectral differentiation and barycentric resampling matrices: our main tools in solving DDEs via spectral collocation. In Section \ref{sec:main} we show how these matrices are combined to solve some simple linear and nonlinear DDEs of various delay types. In Section~\ref{sec:functional} we discuss what small modifications are required so as to apply the method to FDEs. In Section \ref{sec:chebfun} we provide details of how the approach has been incorporated in the Chebfun/Chebop framework for the automated solution of DDEs/FDEs and demonstrate some more challenging examples solved via this software. In Section~\ref{sec:ext} we discuss several extensions, including the use of trigonometric interpolants for PDDEs,  before concluding in Section~\ref{sec:conc}.


\section{Preliminaries}\label{sec:prelim}%
\subsection{Barycentric interpolation formula}\label{subsec:bary} Given 
 data $\{(t_1,y_1),\ldots, (t_n,y_n)\}$ with $\{t_k\}$ finite and distinct, the unique  polynomial interpolant of degree at most $n-1$ may be written in barycentric form as
\begin{equation}\label{eq:barycentric}
p(t) = {\displaystyle\sum_{k=1}^{n}\frac{w_ky_k}{t-t_k}}\bigg/{\displaystyle\sum_{k=1}^{n}\frac{w_k}{t-t_k}},
\end{equation}
where the barycentric weights 
are given by $w_k = C_n\prod_{\ell\neq k}(t_k-t_\ell)^{-1}$ with $C_n$ an arbitrary nonzero constant~\cite[(4.2)]{berrut2004}.\footnote{$C_n$ can be chosen to normalise the weights, for instance, such that $\max_{k}|w_k| = 1$. Constant aside, any other choice of weights will give rise to a rational interpolant which may or may not be free of poles in the interpolation interval.} 
{Stability and conditioning 
of~\eqref{eq:barycentric}
depends on the choice of nodes, $\{t_k\}$, 
which we discuss in Section~\ref{subsec:chebpts}.} For now, let us denote by ${\bm{t}}$ and $\bm{y}$ the $n\times1$ vectors 
$\bm{t} = [t_1, t_2, \ldots, t_n]^\top$ and $\bm{y} = [y_1, y_2, \ldots, y_n]^\top$
. Since \eqref{eq:barycentric} is linear in the data $\{y_k\}$, it can be expressed as $p(t) = P(t;\bm{t})\bm{y}$, where
\begin{equation}\label{eqn:quasimatrix}
 P(t;\bm{t}) = \left[\begin{array}{ccccccc}\ell_1(t) & \ell_2(t) & \cdots & \ell_n(t)\end{array}\right]\ \ \text{and} \ \ \ \ell_k(t) = {\displaystyle\frac{w_k}{t-t_k}}\bigg/{\displaystyle\sum_{k=1}^{n}\displaystyle\frac{w_k}{t-t_k}}.
\end{equation}
Here 
$P(t;\bm{t}) = [\ell_1(t), \dots, \ell_n(t)]$
is a {quasimatrix}, an $\infty$-by-$n$ matrix, whose columns contain the degree $n-1$ Lagrange polynomials for the points $\{t_k\}$. 

By differentiating the columns of $P(t;\bm{t})$ 
and evaluating the result at $t = \bm{t}$, one obtains the pseudospectral differentiation matrix, $D = P'(\bm{t};\bm{t})$ \cite[Section 9.3]{berrut2004}. Code for constructing such matrices is given in Table~\ref{fig:matcode}.
Alternatively, given a different set of points, say $\{\tau_j\}_{j=1}^m$, $P(\bm{\tau};\bm{t})\bm{y}$ evaluates the degree $n$ polynomial interpolant of the data $\{(t_k,y_k)\}_{k=1}^n$ at the $m$ points $\bm{\tau} = [\tau_1,\ldots,\tau_m]^\top$. \cite{Driscoll2016} use this {\em barycentric resampling matrix} with $m < n$ to {downsample} discretised linear operators 
in spectral collocation, leaving space to apply boundary conditions or other side constraints. Here we will use such matrices to evaluate delay terms in DDEs. For instance, a delay term $y(qt-p)$ in a DDE can be discretised as something like $P(q\bm{t}-p;\bm{t})\bm{y}$.\footnote{Here we assume  $p$ and $q$ are such that $q\bm{t}-p$ is contained within the interpolation interval of~(\ref{eq:barycentric}). In general, one  may need to subdivide $[0, T]$ and consider multiple interpolants in order to achieve this. See Section~\ref{subsec:discrete} for details.} Concise MATLAB code for constructing $P(\bm{\tau};\bm{t})$ is given in Table~\ref{fig:matcode}.

\begin{table}[t!]
\begin{tabularx}{\textwidth}{X|l}
\end{tabularx}
\hspace*{35pt}\begin{minipage}[]{.45\textwidth}
\begin{lstlisting}[morekeywords={end}]
function D = Diffmat(t, w)  
    D = w.'./(w.*(t-t.'));
    n = length(t); ii = 1:n+1:n^2;       
    D(ii) = 0;  D(ii) = -sum(D,2); 
end
\end{lstlisting}
\end{minipage}\ \ \ \ \ \ \ 
\begin{minipage}[]{.49\textwidth}
 \begin{lstlisting}[morekeywords={end}]
function P = Barymat(tau, t, w)
    P = w.'./(tau-t.');
    P = P./sum(P, 2);
    P(isnan(P)) = 1;
end
\end{lstlisting} 
\end{minipage}
\caption{\rm  \footnotesize Left:  MATLAB code to generate the pseudospectral differentiation matrix, $D$. Similar code is available in a variety of software libraries, such as {\tt diffmat} in Chebfun~\citep{Chebfun} and {\tt poldiff} in DMSUITE~\citep{DMSUITE}. Right: MATLAB code (modified from \citep{Driscoll2016}) to generate a barycentric resampling matrix, $P(\tau;\bm{t})$. 
}\label{fig:matcode}
\end{table} 



\subsection{Interpolation at Chebyshev points}\label{subsec:chebpts}
It is well known that interpolation at arbitrary nodes may be ill-conditioned. For example, interpolating at uniformly-spaced nodes leads to the Runge phenomenon and exponential ill-conditioning~\citep[Chapter 13]{trefethen2019}. Suitable nodes for polynomial interpolation cluster towards the endpoints of the interpolation interval. 
Of such points, the Chebyshev--Gauss--Lobatto nodes, 
\begin{equation}
t_k = -\cos\left({(k-1)\pi}/(n-1)\right), \quad k = 1, \ldots, n,
 \end{equation}
(which hereafter we simply call  {\em Chebyshev points}) on the interval $[-1,1]$ are the most widely used in spectral methods for nonperiodic problems~\citep{boyd2001, fornberg1996, trefethen2000}, and we focus on these in this manuscript. However, the approach we describe is equally valid in the case of collocation with, say, Gauss--Legendre or Gauss--Jacobi nodes, or even a rational spectral method defined by some alternative choice of the barycentric weights~\citep{berrut2014}. The barycentric weights for polynomial interpolation at the Chebyshev points are known explicitly~\citep{berrut2004}: 
\begin{equation}\label{eqn:chebwts}
w_k = \left\{\begin{array}{llll}
                          1/2, &k = 1,\\
                          (-1)^{k+1}, &k = 2,\ldots,n-1,\\
                          (-1)^{n+1}/2, &k = n.\\
                         \end{array}\right.
\end{equation}

With these points and weights the barycentric formula~\eqref{eq:barycentric} is numerically stable and the interpolation problem is well-conditioned~\citep{higham2004}. If the $\{y_k\}$ are sampled from a function $y(t)$ analytic in a neighbourhood of $[-1,1]$, then the Chebyshev  polynomial interpolant and its derivatives converge geometrically to $y(t)$ as $n$ is increased. In particular, for all $
\nu>0$, $\|\rho^{(\nu)}-y^{(\nu)}\|_{L^\infty_{[-1,1]}} = \mathcal{O}(\rho^{-n})$ as $n\rightarrow
\infty$, where $\rho$ is the radius of the largest Bernstein ellipse with foci at $\pm1$ that lies within the domain of analyticity of $y$~\citep[Theorem 8.2]{trefethen2019}. {Conversely, if $y$ has only $k$ continuous derivatives in $[-1,1]$, then convergence is only algebraic with a rate $\mathcal{O}(n^{-k})$ ~\citep[Theorem 7.2]{trefethen2019}.}

Problems defined on finite intervals other than $[-1,1]$ can be handled by an affine transformation. For instance, if solving on a domain $[0, T]$ then the points $t_k$ are mapped by 
\begin{equation}\label{eqn:cov}
t_k\mapsto \tfrac{T}{2}(t_k+1),
\end{equation}
and the barycentric weights remain unchanged.


\subsection{Spectral collocation for ODEs}
The differentiation matrix, $D$, as computed by the code in Table~\ref{fig:matcode} (left), is the primary tool for spectral collocation~\citep[Chapter 11]{trefethen2000}. Higher-order differentiation matrices can be computed by taking powers of $D$, although there are more accurate and more efficient approaches of calculating higher-order matrices directly~\citep{baltensperger2003}, \cite[(3.28)]{tee2006}. Non-constant coefficients can be incorporated by multiplying by a diagonal matrix containing the values of the coefficient function evaluated at the grid points. Nonlinear problems can be solved by a Newton--Raphson iteration, and piecewise smooth problems by domain decomposition (both of which are discussed below). Boundary conditions can be enforced in a number of ways. The most common approach is boundary bordering (sometimes called row replacement), whereby rows of the discretised operator are removed and replaced with the boundary constraints~\cite[Section 6.4]{boyd2001},~\cite[p.135]{trefethen2000}. An alternative approach is to project the discretised operator onto a lower degree subspace via the barycentric resampling matrix and use the boundary constraints to `square-up' the problem~\citep{Driscoll2016}. For simplicity we employ the more traditional boundary bordering approach in this work, but any reasonable method of enforcing the boundary conditions (including rectangular projection) is equally applicable.

\begin{example}\label{example:1}
For comparison with the DDEs that follow, we briefly recall how a simple ODE can be solved via spectral collocation. In particular, consider the elementary  differential equation
\begin{equation}\label{eqn:example1}
 y'(t) = -y(t), \quad y(0) = 1, 
\end{equation}
with solution $y(t) = e^{-t}$. To solve, the $n\times n$ Chebyshev differentiation matrix  is constructed using the code in Table~\ref{fig:matcode} (left) and the ODE discretised as $A\bm{y} = \bm{b}$ with $A = D+I$. Boundary conditions are enforced by boundary bordering, i.e., removing the first rows of $A$ and $\bm{b}$ and replacing them with the boundary constraint (here the first row of the identity matrix and the scalar 1, respectively). The resulting square system is then solved for $\bm{y}$, the approximate solution on the Chebyshev grid $\bm{t}$. Concise code for solving~\eqref{eqn:example1} on the interval $[0, 1]$ and a plot demonstrating geometric convergence of the discretisation as $n$ is increased are given in Figure~\ref{fig:example1}.  
\begin{figure}[!t]
\hspace*{20pt}\begin{minipage}[]{.4\textwidth}
\begin{lstlisting}
n = 12; y0 = 1;
[t,~,w] = chebpts(n, [0,1]);
z = zeros(n-1, 1);
D = Diffmat(t, w);
I = eye(n); B = I(1,:);
A = D + I;  b = z;  
y = [B;A(2:end,:)]\[y0;b(2:end)]
\end{lstlisting}
\end{minipage}
\begin{minipage}[]{.5\textwidth}\hfill
\includegraphics[height=85pt]{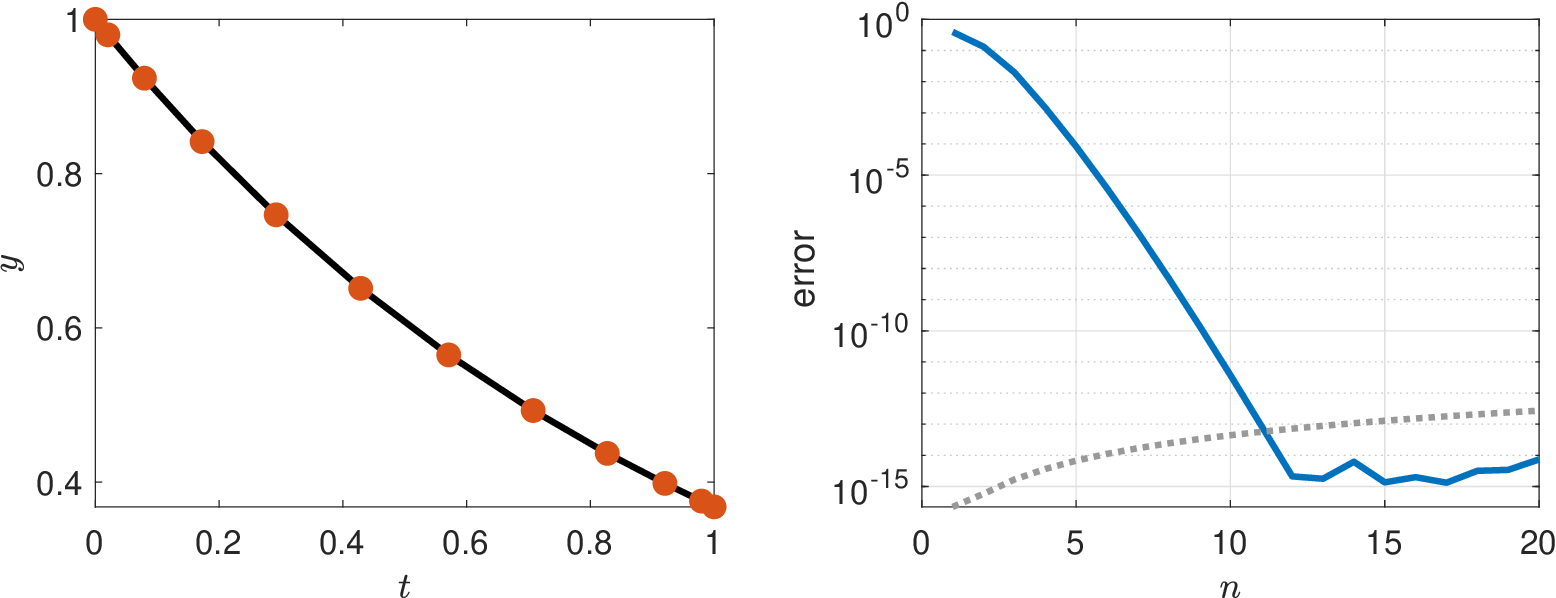}  
\end{minipage}
\caption{Left: MATLAB code for solving the ODE~\eqref{eqn:example1} using Chebyshev spectral collocation. The line {\tt [t,}\lstinline{\~}{\tt,w] = chebpts(n, [0,1])} (from Chebfun) returns the {\tt n}-point Chebyshev grid on the interval $[0,1]$ and the associated barycentric weights (\ref{eqn:chebwts}). Centre: Exact solution (solid) and approximate solution (dots) using an $n =$ 12-point discretisation. Right: The solid line shows geometric convergence of the solution down to around machine precision, $M_\varepsilon = 2^{-52}$,  as $n$ is increased. Here, and throughout, the reported error is the $\ell^\infty$ error at the interpolation points. The dashed line shows (again, here and throughout) $M_\varepsilon\times{\rm cond}(A)$ as a function of $n$. 
}\label{fig:example1}
\end{figure}
\end{example}


\section{Spectral collocation for DDEs}\label{sec:main}%
We now demonstrate how the Chebyshev spectral collocation method for ODEs can be readily extended to solve DDE problems. The examples in this section are intentionally simple, so as to clearly demonstrate the proposed method. In particular, in most cases we choose a specific right hand-side so that the solution is known (typically some form of exponential), allowing us to easily compute the error in the approximation. More interesting and practical examples are given in Section~\ref{sec:chebfun}.

\subsection{Proportional delay}\label{subsec:proportional}
We start by considering proportional (sometimes called pantograph-type) delay of the form $y(qt)$, $q < 1$ on the interval $[0, 1]$.\footnote{Any other finite domain can be dealt with via an affine change of variables as in~(\ref{eqn:cov}). Discrete delays introduced by such a change of variables can be dealt with in the manner described in the next section.} In this setting $qt\in[0, 1]$ and no history function or domain decomposition is required. Such a delay term can therefore be implemented directly by the resampling matrix $P(q\bm{t};\bm{t})$. 
\begin{example}\label{example:2} 
Consider the DDE 
\begin{equation}\label{eqn:example2} 
 y'(t) = - y(t) - y(t/2) + e^{-t/2}, \quad y(0) = 1,
\end{equation} 
with solution $y = e^{-t}$.\footnote{A similar DDE appears in~\citep{ali2009}, where it is solved using a Legendre collocation method. A stiff version of the problem is considered by~\cite{el-safty1998}, where a small parameter is introduced in front of the $y'$ term.} Discretisation of the differential and identity operators proceeds as in Example~\ref{example:1} and  the delay term is discretised using the barycentric resampling matrix as described above. This results now in a system $A\bm{y} = \bm{b}$ where $A = D + I + P(\bm{t}/2;\bm{t})$. 
As before, boundary conditions can be enforced by row replacement.
Code is provided in Figure~\ref{fig:example2} (left) and, as in the ODE case, we observe geometric convergence as the discretisation size is increased (right).  Multiple delays of proportional type, e.g., $y(q_1t)$, $y(q_2t)$, can by included by incorporating additional resampling matrices, $P(q_1\bm{t}; \bm{t})$, $P(q_2\bm{t}; \bm{t})$. Extensions to systems of DDEs and/or nonlinear problems proceed as in the ODE case---see~\citep[Section 4]{Driscoll2016}.

\begin{figure}[!t]
\hspace*{20pt}\begin{minipage}[]{.4\textwidth}
\begin{lstlisting}
n = 12; y0 = 1;
[t,~,w] = chebpts(n, [0,1]);
D = Diffmat(t, w);
I = eye(n); B = I(1,:); 
P = Barymat(t/2, t, w);
A = D + I + P; b = -exp(t/2);
y = [B;A(2:end,:)]\[y0;b(2:end)]; 
\end{lstlisting}
\end{minipage}
\begin{minipage}[]{.5\textwidth}\hfill
\includegraphics[height=85pt]{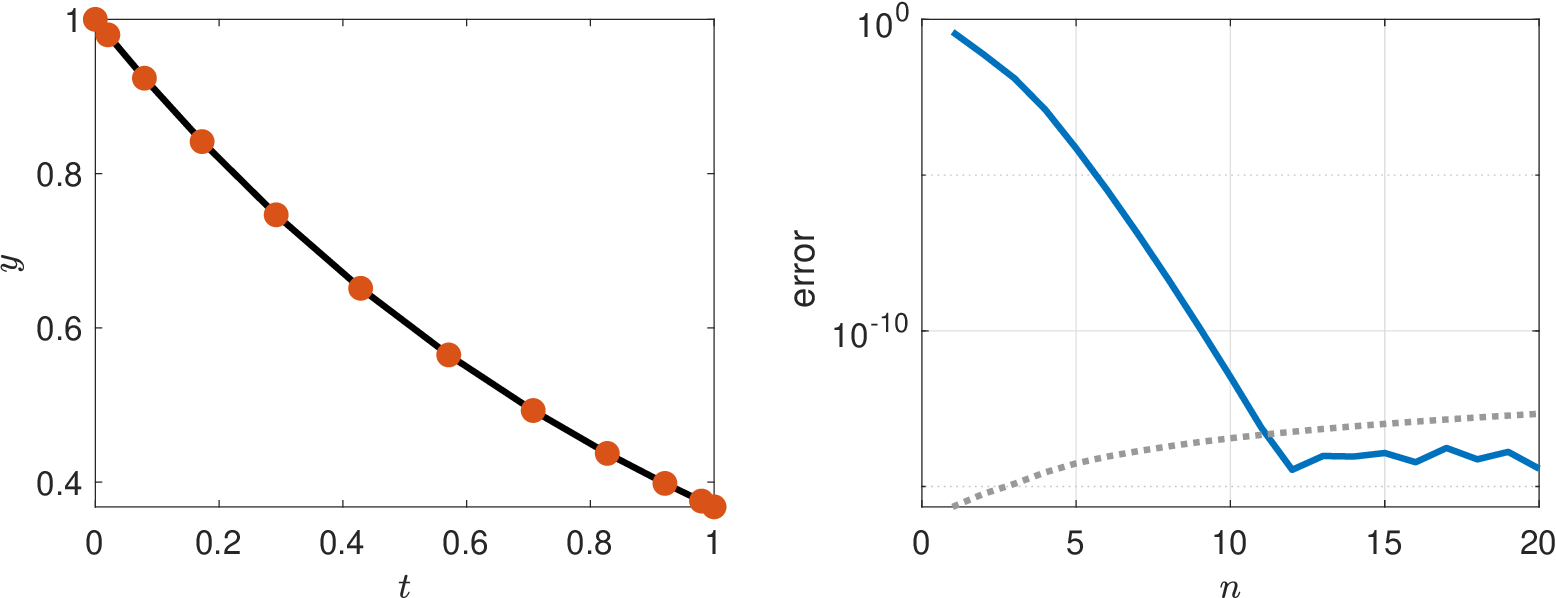} 
\end{minipage}
\caption{Left: MATLAB code for solving the DDE~\eqref{eqn:example2} using Chebyshev spectral collocation. The code is a simple modification of that from Figure~\ref{fig:example1} for solving the ODE~\eqref{eqn:example1}. Centre: Solution using a 12-point Chebyshev grid. Right: Convergence of the solution (solid) and $M_\varepsilon\times{\rm cond}(A)$ (dashed) as $n$ is increased. As in the ODE case, geometric convergence to around the level of machine precision is observed. The condition number of the discretised DDE behaves similarly to the ODE case.}\label{fig:example2} 
\end{figure}

\end{example}

\subsection{Discrete delay and domain decomposition}\label{subsec:discrete}
We next consider discrete delay. For simplicity we start with problems that have only trivial history functions, i.e., $\phi(t) = 0$ for all $t < 0$. Non-trivial history functions can be dealt with in the manner described in the next section.

A challenge arising with discrete delays (or with proportional delays on intervals not including the origin) is that discontinuities caused by non-smooth initial data will propagate into the solution whenever the delay becomes active~\citep{baker1995}, and the rapid convergence of the spectral collocation method will be lost~\citep[Chapter 4]{trefethen2000}. If the delay (or delays)
are short relative to the length of the solution interval, then many discontinuities will be introduced and the use of a spectral method is unlikely to be beneficial. However, if only a few discontinuities are introduced with locations known a priori, then these can be dealt with by a domain decomposition/multidomain approach. In particular, if the solution is piecewise analytic on each subdomain, then geometric convergence may be maintained, as in the case of ODEs~\citep{driscollfornberg1998}. We describe this approach first for a simple problem leading to only one singularity in the domain of interest before outlining the more general approach. A similar approach is employed by~\citep{Borgioli2020} and~\citep{engelborghs2001} in the context of stability analysis of time-periodic delayed dynamical systems (see also Example~\ref{example:11}).

\begin{example}\label{example:3}
Consider the DDE
\begin{equation}\label{eqn:example3}
 y'(t) = -y(t) -y(t-\tfrac12), \quad y(t<0) = 0, \quad y(0) = 1,
\end{equation}
with solution
\begin{equation}
 y = \left\{\begin{array}{llll}
                      e^{-t}, & 0 \le t \le \tfrac12,\\
                      e^{-t+1/2}(\tfrac12 -t + {e^{-1/2}}), & \tfrac12 \le t \le 1.
                     \end{array}\right.
\end{equation} 
From the delay term $y(t - \tfrac12)$, we know a priori that the solution will contain a jump in its first derivative at $t = \tfrac12$ when the delay becomes active (and its second, third, fourth, $\ldots$  derivatives at $t = 1, \tfrac32, 2, \ldots$ if we were to integrate further---see below). As such, we discretise on a piecewise domain with two Chebyshev grids, an $n_L$-point grid,  $\bm{t}_{L}$, on $[0,\tfrac12]$ and an $n_R$-point grid, $\bm{t}_{R}$, on $[\tfrac12,1]$, and seek an approximate solution $\bm{y} = [\bm{y}_L ; \bm{y}_R]$ on each grid. A discretised differential operator is constructed for each subinterval and a continuity condition $[-\bm{e}_{n_L,n_L}^\top, \bm{e}_{n_R,1}^\top]\bm{y} = 0$ (where $
\bm{e}_{n,k}$ denotes the $k$th column of an $n\times n$ identity matrix) enforced so that the solution is continuous at $t = \tfrac12$.\footnote{For higher-order problems, continuity of additional derivatives is enforced at the interface by including appropriate rows of suitable differentiation matrices.} Combined with the initial condition and boundary bordering, this results in a square system that can be solved for $\bm{y}$. Excluding the delay term, details of such domain decomposition are described in detail in~\cite[Section 4.3]{Driscoll2016}.  

To incorporate the delay term, we first note that $y(t - \tfrac12) = 0$ for $t < \tfrac12$, so this term may be neglected in the first subinterval. For any $t$ in the second subinterval, $[\tfrac12,1]$, we note that $t-\tfrac12\in[0, \tfrac12]$; thus the solution on the first subinterval is interpolated to the appropriate points by the resampling matrix $P(\bm{t}_R-\tfrac12, \bm{t}_L)$
.\footnote{If the width and discretisation size of the two intervals is the same, then this will simply be the identity matrix, but in general the resampling matrix will be nontrivial. See Remark~\ref{remark1} for further discussion.} This gives rise to the block linear system in the left panel of Figure~\ref{fig:example3},   where again we use boundary bordering to incorporate the boundary and continuity conditions. For this example, and indeed any IVP-DDE, the resulting system is block lower triangular (BLT) and one can solve first for $\bm{y}_L$ and then for $\bm{y}_R$. In the context of solving DDEs, this is referred to as the {\em method of steps}~\cite[Chapter 1.2]{hale1971}. Whilst this structure could certainly be exploited here, FDEs and boundary-value problem DDES (BVP-DDEs) can be implemented using the same framework as described above, and these will not result in block lower triangular systems. 
Therefore, for simplicity and generality, we solve the whole resulting block discretisation at once using standard direct methods and defer exploration of more efficient domain decomposition strategies for future work. MATLAB code for reproducing the solution to~(\ref{eqn:example3}) (and all other examples in this paper) is available at the online repository~\citep{hale2022}.%
\begin{figure}[!t]
\hspace*{-5pt}\begin{minipage}[]{.46\textwidth}
$$\left[\!\!\!\begin{array}{cccc}
   \bm{e}_{n_L,1}^\top\\[.1em]
   \overline{D_L} + \overline{I} &  \\
   -\bm{e}_{n_L,n_L}^\top & \!\!\! \bm{e}_{n_R,1}^\top\\[.1em]
   \overline{P}(\bm{t}_R-\tfrac12;\bm{t}_L) & \!\!\! \overline{D_R} + \overline{I}
  \end{array}\!\!\!\right]\!
  \left[\!\!\!\begin{array}{cccc}
   \bm{y}_L \\ \bm{y}_R
  \end{array}\!\!\!\right]=
  \left[\!\!\begin{array}{cccc}
   1\\\bm{0}\\0\\\bm{0}  
  \end{array}\!\!\!\right]$$
  \vspace*{1em} 
\end{minipage}
\begin{minipage}[]{.56\textwidth} 
\includegraphics[height=85pt]{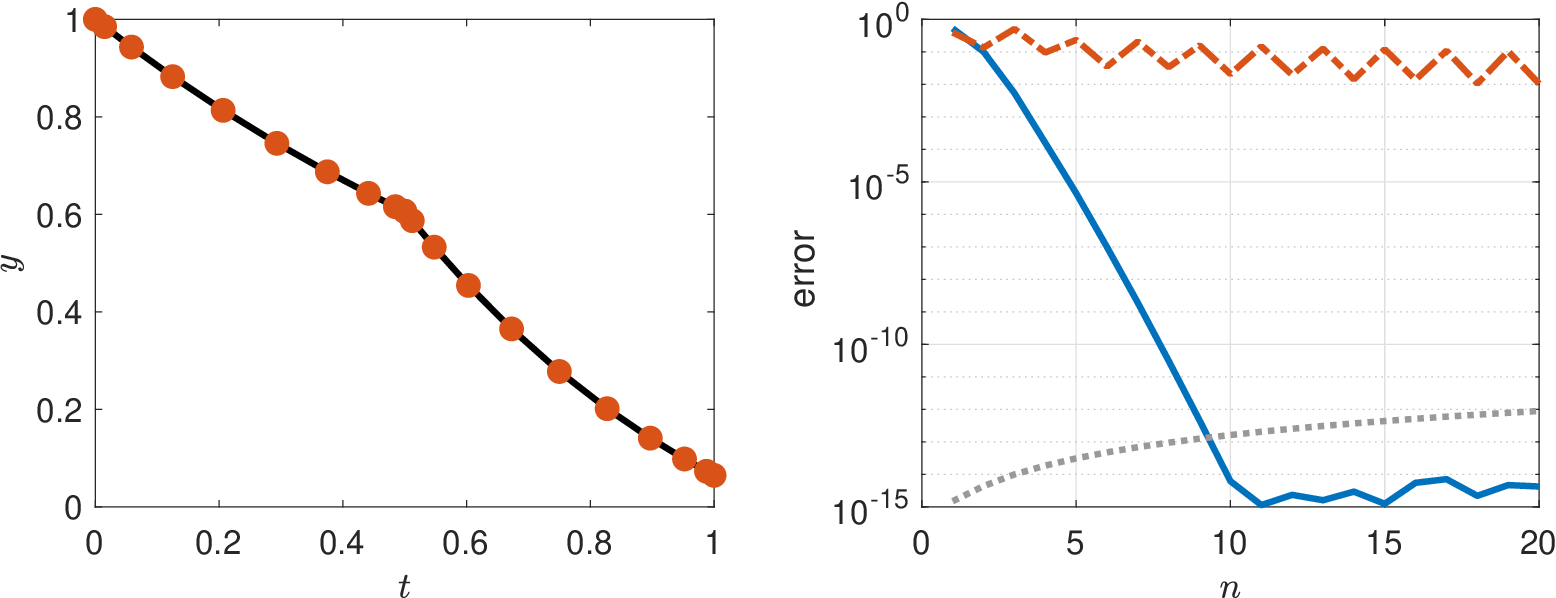}  
\end{minipage}
\caption{Left: Block linear system resulting from the spectral collocation discretisation of~\eqref{eqn:example3}. The overline notation is used to indicate that the first row of the corresponding matrix is removed in order to implement boundary/continuity constraints via boundary bordering. The system is block lower triangular, and one could first solve for $\bm{y}_L$ and then $\bm{y}_R$, akin to the method of steps. However, this will not be the case in context of FDEs as we discuss shortly. Centre: Solution using 10- and 11-point Chebyshev grids on $[0,1/2]$ and $[1/2,1]$, respectively. Right: Convergence as $n$ is increased, where $n$ and $n+1$ are the discretisation sizes of the left and right subintervals, respectively. Using a multidomain discretisation on $[0,1/2]\cup[1/2,1]$ so that the solution is analytic on each subdomain results in geometric convergence (solid line). If only a single Chebyshev grid on $[0, 1]$ is used, then algebraic convergence is observed (dash-dotted line) due to the discontinuity in the derivative of the solution at $t = 1/2$.}\label{fig:example3} 
\end{figure}

The other two panels of Figure~\ref{fig:example3} show the computed piecewise smooth solution $\bm{y} = [\bm{y}_L ; \bm{y}_R]$ (centre) and convergence of the method as the discretisation size is increased (right, solid line), here using an $n$-point discretisation on the left subinterval and an ($n+1$)-point discretisation on the right---see Remark \ref{remark1}. Since the solution is piecewise analytic on the two subintervals, geometric convergence is obtained. Finally, the dashed line in the third figure shows the result of applying a similar approach to discretise this DDE but without introducing the {breakpoint} at $t = \tfrac12$. As expected, the jump in $y'$ at this point results in only algebraic convergence of the numerical solution. As with most numerical methods for DDEs of this type, knowledge of the discontinuity locations is vital~\citep{baker1995}.
\end{example}

\remark{\label{remark1}The choice of using a different discretisation size on each subinterval in the example above is certainly not necessary. We do so to demonstrate that this offers no difficulty to the method. For this example, choosing an $n$-point grid on both subintervals would result in $\overline{P}(\bm{t}_R-\tfrac12;\bm{t}_L)$ simplifying to an identity matrix, thus providing a sparser linear system. However, the ability to use a more general interpolation operator not only allows for differing refinement on each of the subintervals, but is necessary when subintervals are not of a uniform size, which is typically the case when there are two or more discrete delays present or the delay is nonlinear---see Example~\ref{example:5}.}

\vspace*{1em}

Consider now a more general situation and suppose that the propagation of singularities requires that the solution interval $[0, T]$ is decomposed into $m$ subintervals with breakpoints $0 = T_0 < T_1 < \ldots < T_m = T$, and each subinterval $[T_{k-1}, T_k]$ is discretised using an $n_k$-point Chebyshev grid, $\bm{t}_k$, for  $k = 1, \ldots, m$. An evaluation/delay operator $y(\tau(t))$ can then be realised as a block matrix
\begin{equation} P = \left[\begin{array}{cccc}
    P_{1,1} & P_{1,2} & \cdots & P_{1,m}\\
     P_{2,1} & P_{2,2} & \cdots & P_{2,m}\\
     \vdots & \vdots & \ddots & \vdots\\
    P_{m,1} & P_{m,2} & \cdots & P_{m,m}\\
           \end{array}\right],\end{equation}
where each $P_{j,k}$ is an $n_j\times n_k$ matrix defined by
\begin{equation} [P_{j,k}]_{\ell,:} = \left\{\begin{array}{llll}
P(\bm{[\tau}_j]_\ell, \bm{t}_k), \quad & T_k < [\bm{\tau}_j]_\ell \le T_{k+1}, \quad \text{where } \bm{\tau}_j = \tau(\bm{t}_j).\\
\bm{0}^\top, & \rm{otherwise}.
                          \end{array}\right. \end{equation}
Note that, as in Example~\ref{example:3},  if $\tau(t) \le t$ then $P_{j,k} = 0$ for $k > j$, i.e., $P$ is block lower triangular. However, for FDEs this will not be the case in general.
 
\begin{example}\label{example:4}
 We repeat Example 3, but now solve on the domain $[0, 2]$. The propagated singularities appear at $t = 0.5, 1, 1.5$, and so we decompose $[0,2]$ into 4 uniform subintervals. The results are shown in Figure~\ref{fig:4domain}. To demonstrate that one need not use the same discretisation on each subinterval, we have arbitrarily used $n$, $n+1$, $n+2$, and $n+3$-points on each subinterval, respectively---see Remark~\ref{remark1}.
 \begin{figure}[!t]
\hspace*{0pt}\begin{minipage}[]{.46\textwidth}
\small
$$\left[\!\!\!\begin{array}{ccccccccccc}
   {\bm{e_{n_1,1}}}^\top\\
   \overline{D_1} + \overline{I} &  \\
   -{\bm{e}_{n_1,n_1}}^\top & \!\!\! {\bm{e}_{n_2,1}}^\top\\
   \overline{P_{2,1}} & \!\!\! \overline{D_2} + \overline{I}\\
&-{\bm{e}_{n_2,n_2}}^\top &{\bm{e}_{n_3,1}}^\top\\
   &\overline{P_{3,2}} &\overline{D_3} +\overline{I} &  \\
   &&-{\bm{e}_{n_3,n_3}}^\top & \!\!\! {\bm{e}_{n_4,1}}^\top\\
   &&\overline{P_{4,3}} & \!\!\! \overline{D_4} + \overline{I}
  \end{array}\!\!\!\right]\!
  \left[\!\!\!\begin{array}{cccc}
   \bm{y}_1 \\ \bm{y}_2 \\\bm{y}_3 \\ \bm{y}_4
  \end{array}\!\!\!\right]=
  \left[\!\!\begin{array}{cccc}
   1\\\bm{0}\\{0}\\\bm{0}\\0\\\bm{0}\\{0}\\\bm{0}  
  \end{array}\!\!\!\right]$$ 
\end{minipage} 
\begin{minipage}[]{.56\textwidth}  
\includegraphics[height=90pt,trim={1.5cm 0 0 -10pt},clip]{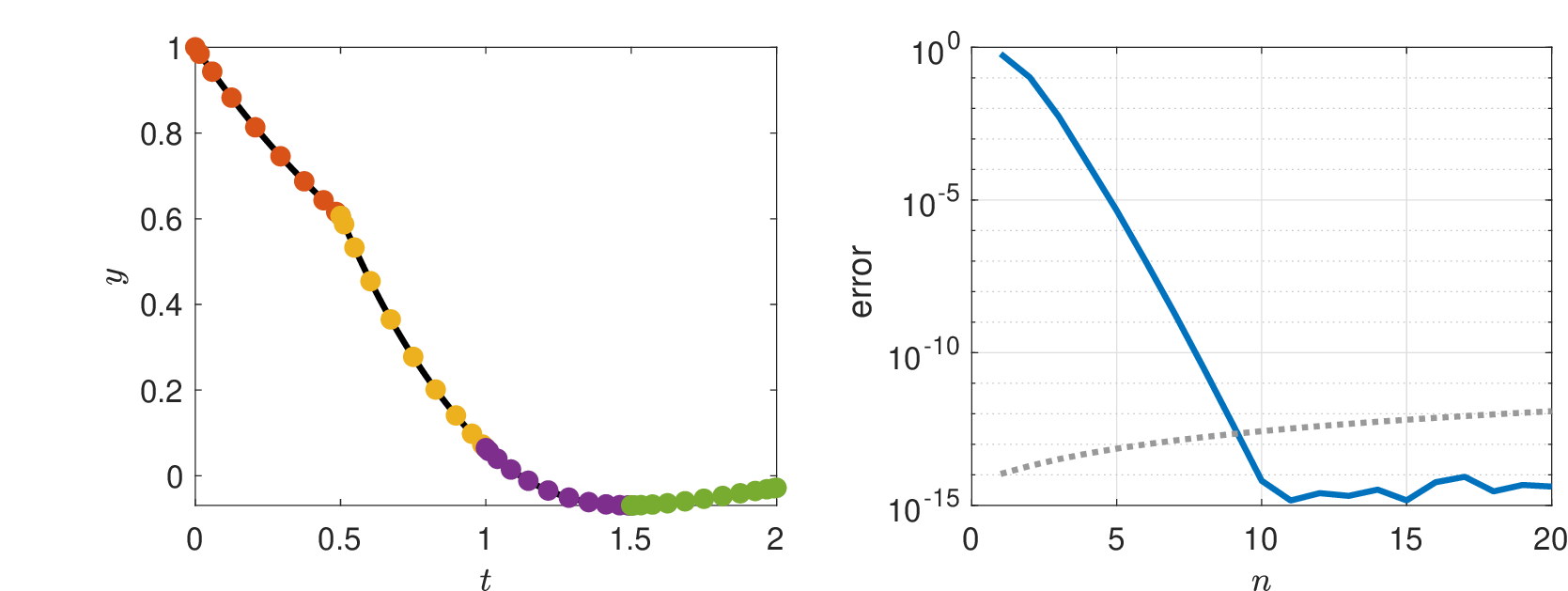}   
\end{minipage}
\caption{Left: Block discretisation of the DDE~(\ref{eqn:example3}) on the interval $[0, 2]$ with `breakpoints' at $0.5, 1, 1.5$.  Here $P_{j,k}$ is given by $P(\bm{t}_j-\tfrac12;\bm{t}_k)$ and ${\bm{e}_n}_k$ denotes the $j$th column of the $n\times n$ identity matrix. Each nonzero entry will result in an essentially dense $n\times n$ block. Centre: Solution using 10-, 11-, 12-, 13-point grids on the respective subintervals (see Remark~\ref{remark1}). Right: Geometric convergence as discretisation size is increased ($n$, $n+1$, $n+2$, $n+3$-point discretisations on each subinterval, respectively.)}\label{fig:4domain}
\end{figure}

Although the solution may look smooth at $t = 1$ and $t = 1.5$, there are discontinuities in the 2nd and 3rd derivatives of the solution at these points, respectively.  However, since the solution is piecewise analytic between these points, we again observe geometric convergence in the numerical solution when using four subintervals with the discontinuity locations as breakpoints. Note again that the system is block lower triangular, so one can emulate the method of steps in solving first for $\bm{y}_1$, then $\bm{y}_2$, and so on.

\remark For many DDEs, the order of jump discontinuities are typically smoothed out along the integration interval. 
For example, in the DDE above, one can readily show that number of continuous derivatives increases by one at each break point, $T_k = k/2$. In principle, for large enough $k$, one could ignore the singularities and sacrifice geometric convergence using small subintervals for high-order algebraic convergence on a larger interval. However, although this phenomenon is fairly common in standard DDE problems with discrete delay, the situation is more complicated for, say, neutral delay-type DDEs~\citep{baker2006}. To avoid undue over-complication, we present this approach as an option, but not pursue it further.  
\end{example}

\subsection{Non-trivial history functions} In many applications the DDE history function, $\phi(t)$, is constant and consistent with the initial condition so that $y(t\le 0) = y_0$. In this case, one can simply consider an equivalent delay $\tau(t)\mapsto \max(\tau(t), 0)$ and enforce a constant history function in a similar manner as described in the previous section. More generally, at least for state-independent delays, history functions can be dealt with by rewriting the DE in a piecewise-defined form with a constant history function. In particular, a DDE in the form of \eqref{eqn:DDE}  can be rewritten as
\begin{equation}
 y'(t)  = \left\{\begin{array}{lll}
            f(t, y(t), \phi(t)), \quad &\tau(t) < 0, \\
            f(t, y(t), y(\tau(t)), &0 \le \tau(t), \\
                \end{array}\right.
                \end{equation} 
with $y(0) = y_0$. Multiple and neutral-type delays can be dealt with in a similar manner. To ensure exponential convergence, a multidomain approach with breakpoints included whenever these delays become active can be employed, as described in the previous section.

\begin{remark}
An important special case is periodic DDEs (PDDEs), which  arise in the context of determining limit cycles of delayed dynamical systems~\citep{michiels2007}. The approach above can be modified to solve such problems by considering $\phi(t) = y({\rm mod}(\tau(t),T))$, where $T$ is the period. This is similar to the methods of~\citep{butcher2011},~\citep{Borgioli2020}, and~\citep{engelborghs2001}, who use (typically low-order) polynomial-based collocation methods to solve PDDEs. In most cases, the period $T$ is then unknown, and must be solved for as part of the problem. However, the periodicity assures that the initial data will be smooth, meaning that a multidomain approach is not typically required in situations. We consider such a situation in Examples~\ref{example:11} and~\ref{example:13} below.
\end{remark}
\subsection{Nonlinear delay}
The approaches described above can be applied to  arbitrary smooth delays functions, $\tau(t) \le t$. In such instances, locating the propagated discontinuities can be challenging, but in principle it can be done~\citep{neves1976}. With these locations known, the domain can be subdivided and the multidomain approach used as described above. 
\begin{example}\label{example:5}
We demonstrate on a simple example with the nonlinear delay $\tau(t) = t^2 - 1/4$. In particular, consider the DDE 
\begin{equation}\label{eqn:historyexample}
 y'(t) = -y(t) - y(t^2-\tfrac14), \quad y(t<0) = 0,  \quad y(0) = 1.
\end{equation}
Here, as in the previous example, a discontinuity in the derivative of the solution is introduced at $t=1/2$, where $\tau(t) = 0$. There will then be an additional discontinuity introduced in the second derivative where $\tau(t) = 1/2$, i.e., at $t = \sqrt{3}/2$. In this manner, all discontinuities can be traced out, even for non-constant delays. For this example, if we solve for $0 \le t \le 1$, these are the only two discontinuities required. (The next occurs when $\tau(t) = \sqrt{3}/2$, i.e., at $t = \tfrac12\sqrt{2\sqrt{3}+1}\approx 1.05$.) Since the subintervals $[0, 1/2]$, $[1/2, \sqrt{3}/2]$, and $[\sqrt{3}/2, 1]$ are of differing widths, even if the same discretisation size is used on each, the resampling matrices are not trivial. We therefore discretise using the approach outlined in Section~\ref{subsec:discrete}. The left panel of Figure~\ref{fig:example5} shows the result of solving with an $n = 12$-point discretisation on each subinterval. The centre panel demonstrates geometric convergence as $n$ is increased.\footnote{Using the method of steps, an explicit (but lengthy) solution to~(\ref{eqn:historyexample}) can be expressed in terms of the integral of an error function. We omit the details.} The final panel shows the sparsity plot of the $n=12$ discretisation.

\begin{figure}[!t]
\centering 
\hspace*{-28pt}\includegraphics[height=105pt,trim={0cm 0 27.5cm 0},clip]{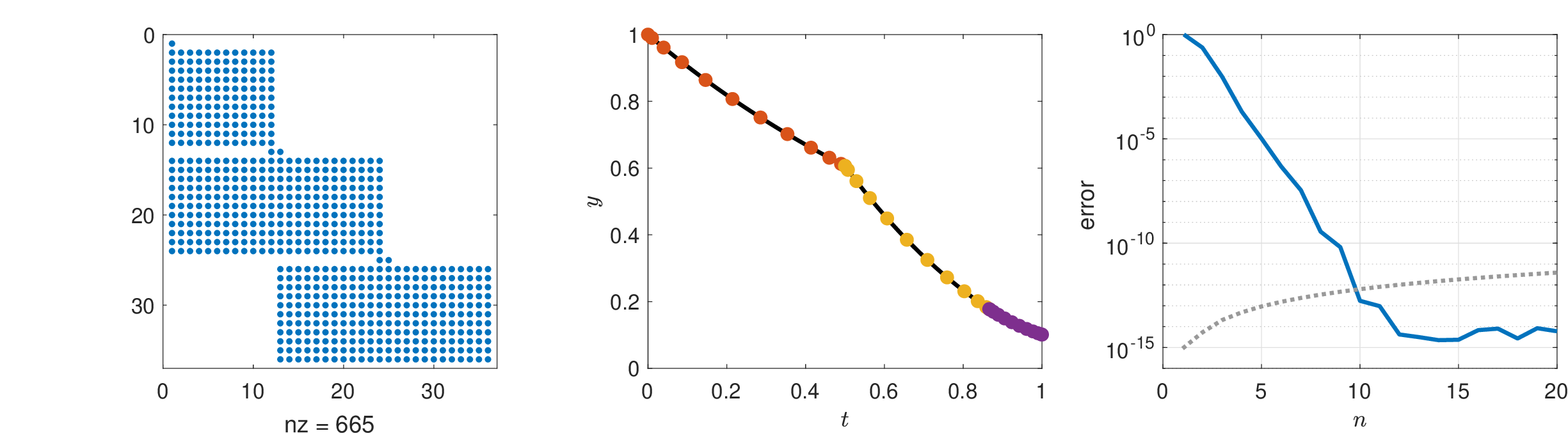}
\includegraphics[height=105pt,trim={15cm 0 0 0},clip]{figures/example5}
\caption{Left: Sparsity plot of the 12-point Chebyshev spectral discretisation. Here, even though the same value of $n$ is used on each subinterval, the interpolation operator is non-trivial (because the intervals are of differing widths), resulting in dense off-diagonal blocks of the discretisation.  Centre: Analytic and computed solution to~(\ref{eqn:historyexample}) on a 12-point grid. Right: Geometric convergence as the discretisation size is increased.}\label{fig:example5} 
\end{figure}  
\end{example}%

\subsection{Nonlinear problems \& state-dependent delay}
In spectral methods for ODEs, nonlinear problems are typically solved by a Newton--Raphson iteration, or some variant thereof---see~\citep[Section 4.2]{Driscoll2016} for details. DDEs allow for an additional form of nonlinearity, whereby the nonlinearity occurs in the delay. This is  known as {\em state-dependent} delay~\citep{Hartung2006}.

Suppose a term of the form $y(g(y))$ appears in the DDE, where $g$ is some differentiable function defined on the range of $y$. This can then be discretised using the barycentric resampling matrix on a grid $\bm{t}$  as $P(g(\bm{y});\bm{t})\bm{y}$, where $\bm{y} = y(\bm{t})$. 
We must linearise this expression about $\bm{y}$ to obtain the Jacobian.

To proceed, we note that, by definition, 
\begin{eqnarray}
	[\nabla \big(P(g(\bm{y});\bm{t})\bm{y}\big)]_{i,j} 
	& = & \frac{\partial}{\partial y_j} [\nabla \big(P(g(\bm{y});\bm{t})\bm{y}\big)]_i\\
	& = &  \frac{\partial}{\partial y_j}\sum_{k=1}^n [P(g(
	\bm{y});\bm{t})]_{i,k} y_k.
\end{eqnarray}	
Applying the product rule, we find
\begin{eqnarray}
	[\nabla \big(P(g(\bm{y});\bm{t})\bm{y}\big)]_{i,j}  &=& \sum_{k=1}^n \Biggl\{\frac{\partial}{\partial y_j}[P(g(
	\bm{y});\bm{t})]_{i,k} y_k + [P(g(\bm{y});\bm{t})]_{i,k} \frac{\partial}{\partial y_j}y_k\Biggr\}\\
	&=& \sum_{k=1}^n \Biggl\{\Biggl[\frac{\partial}{\partial y_j}\frac{\frac{w_k}{g(y_i) - t_k}}{\sum_{\ell=1}^n \frac{w_\ell}{g(y_i) - t_\ell}}\Biggr] y_k\Biggr\}+[P(g(    	\bm{y});\bm{t})]_{i,j}\\
		&=& \delta_{i,j}\Biggl\{\Biggl[\frac{d}{d y_i}\sum_{k=1}^n \frac{\frac{w_k}{g(y_i) - t_k}}{\sum_{\ell=1}^n \frac{w_\ell}{g(y_i) - t_\ell}}\Biggr] y_k\Biggr\}+[P(g(    	\bm{y});\bm{t})]_{i,j},\label{eqn:curly}
\end{eqnarray}
where $\delta_{i,j}$ is the Kronecker delta.
Comparing the term in braces in~(\ref{eqn:curly}) to~(\ref{eq:barycentric}), we see that the former is precisely the derivative of the polynomial interpolant to $\{t_k, y_k\}$, evaluated at $g(y_i)$, hence
\begin{eqnarray}
[\nabla \big(P(g(\bm{y});\bm{t})\bm{y}\big)]_{i,j} 	= \delta_{i,j}\frac{d}{dy_i}p(g(y_i))+[P(g(    	\bm{y});\bm{t})]_{i,j}.
\end{eqnarray}
Finally, by the chain rule, we have
\begin{eqnarray}
[\nabla \big(P(g(\bm{y});\bm{t})\bm{y}\big)]_{i,j}
	& = & \delta_{i,j}g'(y_i)p'(g(y_i))+[P(g(    	\bm{y});\bm{t})]_{i,j}\\
	& = &  \delta_{i,j}g'(y_i)P'(g(y_i);\bm{t})\bm{y}+[P(g(    	\bm{y});\bm{t})]_{i,j},
\end{eqnarray}
or, in matrix form, \begin{equation}\label{eqn:linearise0}
\nabla \big(P(g(\bm{y});\bm{t})\bm{y}\big)  =   
 \text{diag}\Big(g'(\bm{y})\cdot\big(P'(g(\bm{y});\bm{t})\bm{y}\big)\Big) + P(g(\bm{y});\bm{t}),
\end{equation}
where $\cdot$ denotes pointwise multiplication.
Rather than construct $P'(g(\bm{y});\bm{t})$ explicitly, we observe that $P'(g(\bm{y});\bm{t}) = P(g(\bm{y});\bm{t})D$, and therefore that
\begin{equation}\label{eqn:linearise}
		\nabla \big(P(g(\bm{y});\bm{t})\bm{y}\big)  =   
		 \text{diag}\Big(g'(\bm{y})\cdot\big(P(g(\bm{y});\bm{t})D\bm{y}\big)\Big) +  P(g(\bm{y});\bm{t}),
\end{equation} 
which can then be incorporated in a Newton--Raphson iteration to solve state-dependent problems. 

A detailed description of the existence and uniqueness of solutions to state-dependent DDEs is beyond the scope of this paper. Furthermore, the propagation of discontinuities is significantly more complex than in the constant delay case~\citep{neves1976}. We put aside these challenges for now and instead consider below a reasonably simple case for which a known (smooth) solution exists. A less contrived example of state-dependent delay is demonstrated in Section~\ref{sec:chebfun}. 

\begin{example}\label{example:6} Consider the DDE
\begin{equation}\label{eqn:example6}
 y'(t) = -y(y(t)) + f(t), \quad t\in[0, 1], \qquad y(0) = 0.
\end{equation}
We employ the { method of manufactured solutions}, taking $f(t) = \cos(t)+\sin(\sin(t))$ so that $y(t) = \sin(t)$. Figure~\ref{fig:example6} demonstrates the result of applying the approach outlined above to the DDE~(\ref{eqn:example6}), using an initial guess $y(t) = t$ in the Newton--Raphson iteration. The table on the left shows  quadratic convergence of the iterations in both the norm of the residual and the norm of the Newton update for an $n = 12$-point discretisation. As with previous examples, the centre panel shows the resulting computed solution and the rightmost panel demonstrates geometric convergence as $n$ is increased.
\begin{figure}[!t]
\hspace*{10pt}\begin{minipage}[]{.4\textwidth}
%
\begin{tabular}{ccc}
	  $k$ &  Residual norm &  Update norm\\\hline
	  0 &  0.71407355247 &  0.26232516612\\
	  1 &  0.05480002458 &  0.01314905164\\
	  2 &  0.00016794991 &  0.00002292528\\
	  3 &  0.00000000051 &  0.00000000004\\[.5em]
\end{tabular}
\end{minipage}\hspace*{12pt}\hfill
\begin{minipage}[]{.56\textwidth}
\includegraphics[height=85pt]{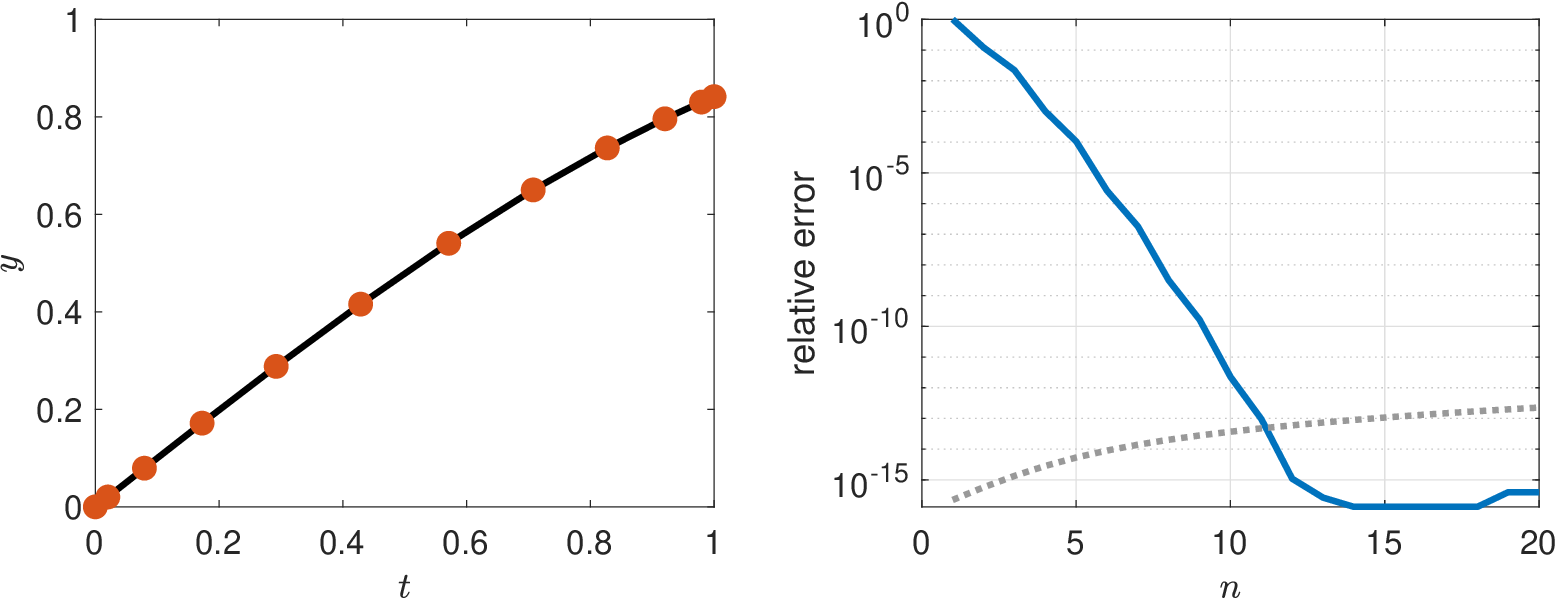}   
\end{minipage}
\caption{Computed solution to the state-dependent delay problem,~(\ref{eqn:example6}). Left: Quadratic convergence of the Newton iteration applied to the $n = 12$ discretisation with initial guess $y = t$. Here $k$ is the iteration number and the norms are the $\ell^\infty$ norms on the grid points. Centre: Computed solution taking $n = 12$. Right: Geometric convergence as $n$ is increased (solid) and $M_\varepsilon$ times the condition number of the Jacobian in the final iteration of the Newton iteration (dashed).}\label{fig:example6}  
\end{figure}  
\end{example} 

\begin{remark}
 A careful reader may have noticed that in the above example we carefully chose $y(t) = \sin(t)$ so that the state-dependent delay satisfied $\tau(t, y(t)) \le t$, and was indeed a {\em delay} differential equation. However, as with previous examples, this is not a requirement of the method, and a state-dependent problem with $\tau(t, y(t)) \not< t$ is considered in Section~\ref{sec:functional}. 
\end{remark}


\subsection{Continuous delay}
Continuous  delays (sometimes called {distributed} delays) result in Volterra-type integro-differential equations. \cite{driscoll2010} presents a detailed description of a Chebyshev spectral collocation method for such equations, and we refer the reader to that work for further information. We include an example here both for completeness and  to demonstrate that, since both methods are natural extensions to standard Chebyshev spectral methods for ODEs, they can be effortlessly combined to solve problems with a combination of both discrete, proportional, and continuous delay types.

\begin{example}\label{example:7}
Consider the DDE below, which combines both proportional and continuous delay:
\begin{equation}\label{eqn:example7}
 y'(t) + \tfrac12y(t/2) = \int_0^t y(s)e^{-(t-s)^2}\,ds, \qquad y(0) = 1.
\end{equation}
The code in the left panel of Figure~\ref{fig:example7} demonstrates compact MATLAB code for solving this problem, combining the resampling approach of Section~\ref{subsec:proportional} with code from \citep{driscoll2010} for discretising the Volterra operator. The computed solution is depicted in the centre panel and geometric convergence as the discretisation size is increased is demonstrated on the right. 
\begin{figure}[!t]
\hspace*{20pt}\begin{minipage}[]{.4\textwidth}
\begin{lstlisting}
K = @(s,t) exp(-(t-s).^2);
n = 14; dom = [0 1]; 
t = chebpts(n, dom);  
D = diffmat(n, dom);            
P = Barymat(t/2, t);
V = K(t, t').*cumsummat(n, dom);  
A = D + P/2 - V ;  rhs = 0*t;     
A(1,:) = eye(1,n); rhs(1) = 1; 
y = A\rhs;
\end{lstlisting}
\end{minipage}\hfill
\begin{minipage}[]{.56\textwidth}
\includegraphics[height=85pt]{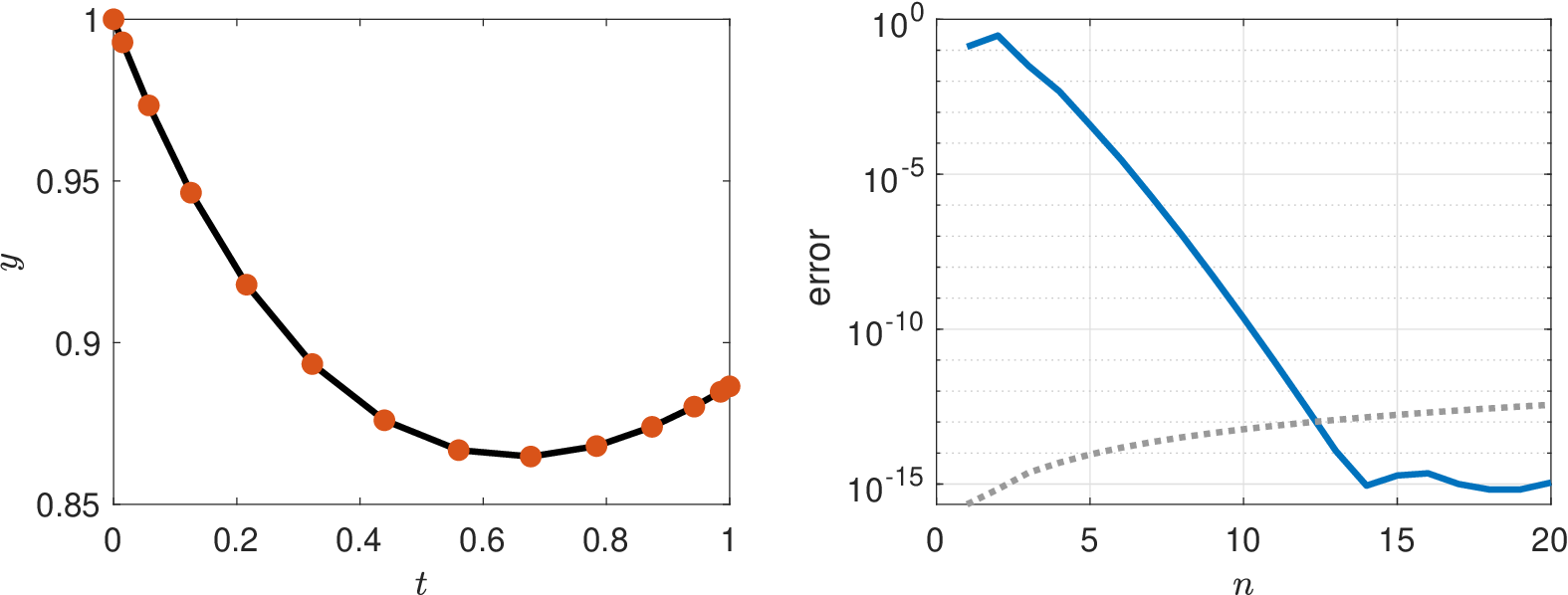}   
\end{minipage}
\caption{Left: MATLAB code for solving the combined proportional and continuous DDE~(\ref{eqn:example7}). The Chebfun function \lstinline{cumsummat(n,dom)}
computes an \lstinline{n}-by-\lstinline{n} spectral collocation matrix representing indefinite integration. Note that here we use the Chebfun implementations of \lstinline{diffmat} and \lstinline{barymat}, which are hardcoded for Chebyshev points. Centre: Computed solution on an 14-point grid. Right: Geometric convergence as the discretisation size, $n$, is increased. Here the $n=30$ discretisation is used a proxy for the unknown true solution.}\label{fig:example7} 
\end{figure}  
\end{example}%
\subsection{Other types of DDE problems}
We have described above how the barycentric resampling approach can be incorporated into standard spectral collocation methods to solve the most common types of DDE problems. The approach can be further generalised to neutral-type DDEs, with a discrete operator $P(\bm{\tau}, \bm{t})D$ representing $y'(\tau(t))$, as well as higher-order problems, DDE-BVPs, systems of DDEs, DDE-EVP, and partial delay differential equations using natural extensions of the techniques described in~\cite[Sections 4--6.]{Driscoll2016} for ODEs. Examples of some of these are given in Section~\ref{sec:chebfun} below.
 
\section{Spectral collocation for functional differential equations}\label{sec:functional}
Unlike time-marching methods, the approach introduced in Section~\ref{sec:main} did not require $\tau(t) \le t$, making it applicable to more general functional differential equations (FDEs). As with DDEs, one must account for propagation of discontinuities caused by non-smooth initial data, as well as the potential need for an additional `history' function for $t > T$~\citep{hale1971}. For the sake of brevity we avoid these complications here, which can be treated in a manner similar to those described in the previous section, and instead consider two examples for which $\tau(t)\in[0, T]$.

\begin{example}
Consider the following FDE:
 \begin{equation}\label{eqn:example8}
  y'(t) = -y(t) - y(1-t^2) + f(t), \quad t\in[0, 1], \qquad y(0) = 1,
 \end{equation}
where here $f(t)$ is chosen as $e^{t^2-1}$ so that~(\ref{eqn:example8}) has the exact solution $y(t) = e^ {-t}$. Clearly  the $1-t^2$ in the argument of the second term on the right hand-side does not satisfy $1-t^2 < t$ for all $0 \le t \le 1$, making this an FDE but not a DDE. However, as shown in Figure~\ref{fig:example8}, discretisation proceeds almost identically to the DDE~(\ref{eqn:example2}) and geometric convergence is again observed as the discretisation size is increased.
\begin{figure}[!t]
\hspace*{20pt}\begin{minipage}[]{.4\textwidth}
\begin{lstlisting}[deletekeywords={end}]
n = 12; y0 = 1;
[t,~,w] = chebpts(n, [0,1]);
D = Diffmat(t, w);
I = eye(n); B = I(1,:); 
P = Barymat(1-t.^2, t, w);
A = D + I + P; b = exp(t.^2-1);
u = [B;A(2:end,:)]\[y0;b(2:end)]; 
\end{lstlisting}
\vspace*{.5em}
\end{minipage}
\begin{minipage}[]{.59\textwidth}
\includegraphics[height=85pt]{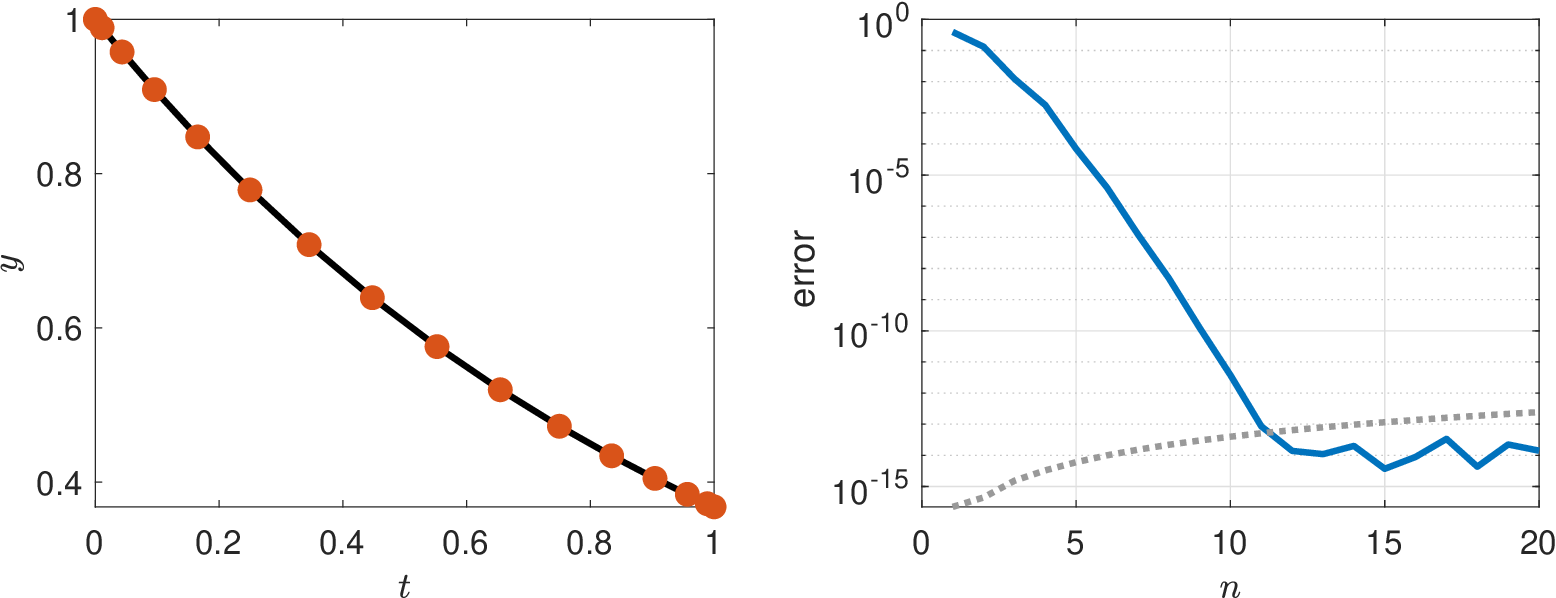}  
\end{minipage}
\caption{Left: MATLAB code for solving the FDE~\eqref{eqn:example8} using Chebyshev spectral collocation. The only change from solving the DDE~(\ref{eqn:example2}) is the first argument of \lstinline{Barymat} on line 5 and the definition of the right-hand side vector, \lstinline{b}. Centre: Solution using 12-point Chebyshev grid. Right: Convergence as $n$ is increased. As in the DDE case, geometric convergence to around the level of machine precision is observed.}\label{fig:example8} 
\end{figure}

\end{example}

\begin{example}\label{example:9}
Here we consider an FDE variant of the DDE in Example~\ref{example:6} in which $y(t) \not< t$, namely\footnote{This FDE was the original motivation for this work after a student in the author's undergraduate class became confused when solving the much more simple ODE $y' = -y^2$, $y(0) =1 $ and asked how one might solve such an equation as~(\ref{eqn:example9}).}\footnote{A DDE closely related to~(\ref{eqn:example9}) is considered in some detail by~\cite{eder1984}.}
 \begin{equation}\label{eqn:example9}
  y'(t) = -y(y(t)), \qquad y(0) = 1.
 \end{equation}
Here we are lucky in that $y(t)\in[0,1]$ for $t\in[0,1]$, meaning that no history functions or domain subdivisions are required. Again, the fact that this is an FDE rather than a DDE has almost no impact on the discretisation, and the implementation proceeds almost identically to that of Example~\ref{example:6}. Figure~\ref{fig:example9} shows convergence of the Newton iteration, the computed solution, and its convergence.%
\begin{figure}[!t]
\hspace*{8pt}\begin{minipage}[]{.4\textwidth}
\begin{tabular}{ccc}
	$k$ &  Residual norm &  Update norm\\\hline
     0  & 1.00000000000 &  1.075290658380\\
     1 &  0.25000000000 &  0.159726357356\\
     2 &  0.00686128071 &  0.002791677486\\
     3 &  0.00000843021 &  0.000005995919\\
     4 &  0.00000000002 &  0.000000000006\\[1em]
\end{tabular}
\end{minipage}\hspace*{14pt}\hfill
\begin{minipage}[]{.59\textwidth}
\includegraphics[height=85pt]{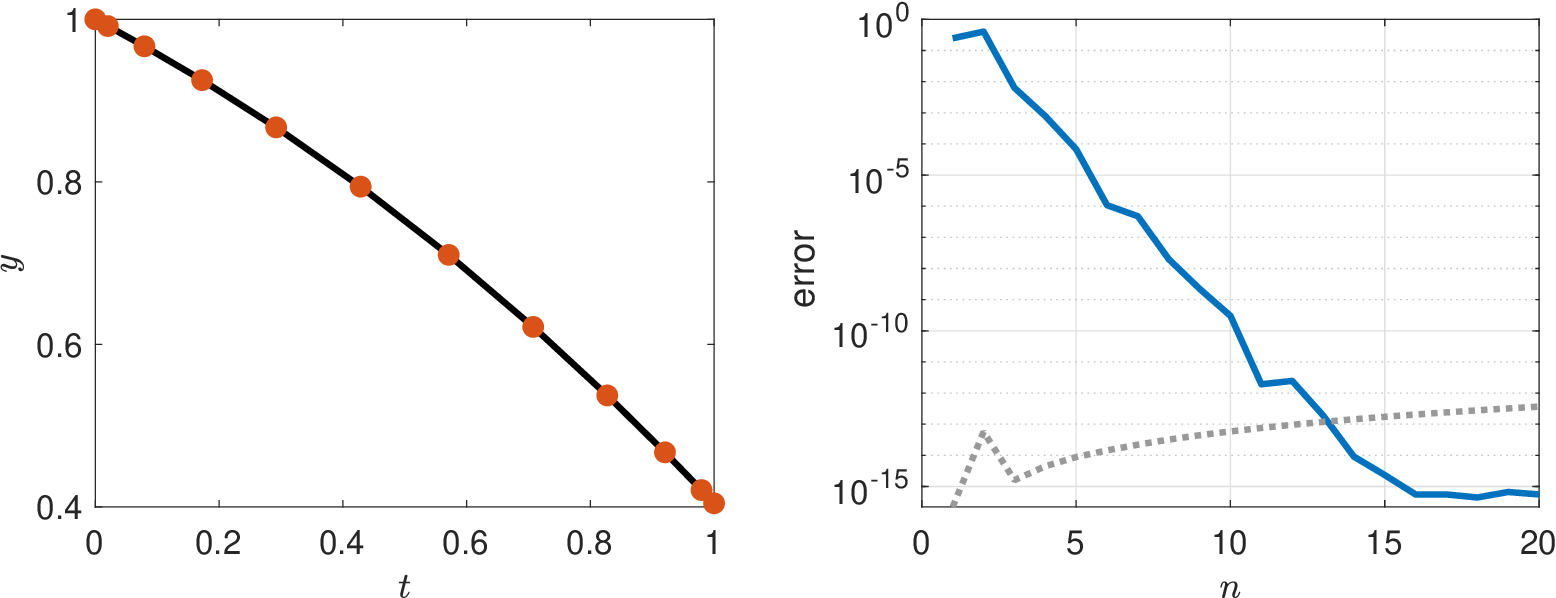}   
\end{minipage}
\caption{Computed solution to the state-dependent FDE~(\ref{eqn:example9}). Left: Quadratic convergence of the Newton iteration applied to the $n = 12$ discretisation with initial guess $y = 1$. Centre: Computed solution taking $n = 12$. Right: Geometric convergence as $n$ is increased, taking the $n = 30$ solution as a proxy for the unknown true solution. }\label{fig:example9} 
\end{figure}  
\end{example} 


\section{Chebfun \& more examples}\label{sec:chebfun}%
Chebfun is an open-source software system for numerical computation with functions~\citep{Chebfun}. Its {Chebop} class and associated methods allow Chebfun to solve a variety of linear and nonlinear ODE boundary-value problems through the use of adaptive spectral  collocation, automatic differentiation, and robust Newton--Raphson iteration~\citep{Driscoll2008, Birkisson2012}. Chebfun can also solve integral equations~\citep{driscoll2010}, differential eigenvalue problems~\citep{Driscoll2008}, as well as initial-value PDEs through the use of method-of-lines type discretisations~\citep[Chapter 22]{exploringodes}.      

The framework outlined in this manuscript uses the same underlying Chebyshev spectral collocation discretisation as Chebop uses for ODEs, allowing the natural extension of Chebop to solve linear FDEs. Furthermore, Chebop's automatic differentiation capabilities are easily modified to allow linearisation of state-dependant delays using~(\ref{eqn:linearise}), thus also allowing the solution of nonlinear FDEs. A detailed description of the software engineering involved in this extension is beyond the scope of this work, but the syntax, usability, and versatility of using Chebfun/Chebop to solve DDEs is demonstrated in various examples below. In addition, Chebfun implementations of the examples from Section~\ref{sec:main} can be found at the repository~\citep{hale2022}.

\subsection{\textbf{\em Chebfun example 1.} (proportional and continuous delay)}
\begin{figure}[!t] 
\hspace*{20pt}\begin{minipage}[]{.65\textwidth}
\begin{lstlisting}
dom = [0 20]; q = 1/2;
N = chebop(@(t,y) diff(y) - 1/100*(q*t-t-10)*y(q*t) - ...
     1/100*(t+20)*exp(-1) - 1/100*cumsum(y) - ...
     1/1000*feval(volt(@(x,s)(x/q-s), y), q*t), dom);
N.lbc = exp(-1);
y = N\0;
sol = chebfun(@(t) exp(t/10-1), dom); 
err = norm(y-sol, inf)./norm(sol, inf) 
\end{lstlisting}
\end{minipage}\hspace*{15pt}
\begin{minipage}[]{.19\textwidth}\vspace*{7pt}
 \includegraphics[height=85pt]{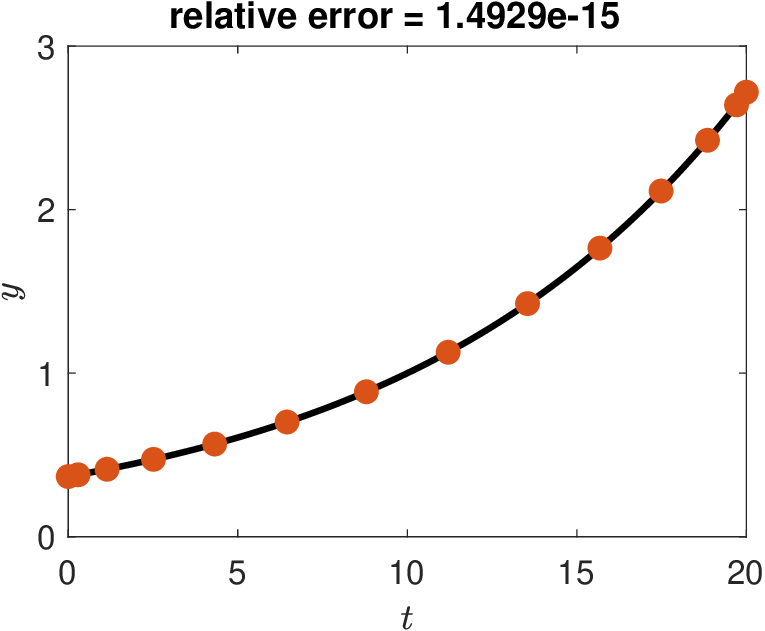} 
\end{minipage}
\caption{Left: Chebfun code to solve the continuous-delay DDE~(\ref{eqn:chebxample1}). Right: Computed solution (dots) and relative error.}
\label{fig:chebxample1}
\end{figure}

Our first Chebfun example is a linear DDE with both proportional and continuous delay, namely
\begin{equation}\label{eqn:chebxample1}
 y'(t) = \frac{1}{100}(qt-t-10)y(qt)+\frac{1}{100}(t+20)e^{-1}+\frac{1}{100}\int_0^ty(s)\,ds + \frac{1}{1000}\int_0^{qt}(t-\tau)y(\tau)\,d\tau,
\end{equation} 
for $t\in[0, 20]$ with $y(0) = e^{-1}$, with exact solution $y(t) = e^{\frac{1}{10}t-1}$.  This example appears in~\citep{wei2012legendre}, where it is solved using a Legendre-based spectral collocation method. The short code in Figure~\ref{fig:chebxample1} demonstrates the ease with which this DDE can be solved in Chebfun. Here the computed solution, taking $q = \frac{1}{2}$, is returned as a length 14 chebfun and has a relative error of around $1.5\times10^{-15}$.

\subsection{\textbf{\em Chebfun example 2.} (nonlinear, neutral type)}\label{sec:chebfunexample2}  
\begin{figure}[!t]  
\hspace*{30pt}
\begin{minipage}[]{\textwidth}
	\begin{lstlisting}[morekeywords={end}]
N = chebop(@(t,y) diff(y) - (2*cos(2*t)*y(t/2)^(2*cos(t)) + ...
                  log(deriv(y,t/2)) - log(2*cos(t)) - sin(t)), [0,.1]);
N.lbc = 1; 
for s = [2, -lambertw(-2/exp(2))]
    N.init = chebfun(@(t) 1+s*t, [0,.1]);
    u = N\0;
    figure, plot(u), title(['s = ' num2str(s, 14)]) 
end 
\end{lstlisting} 
\end{minipage}
\begin{minipage}[]{\textwidth}
\centering
\includegraphics[height=125pt]{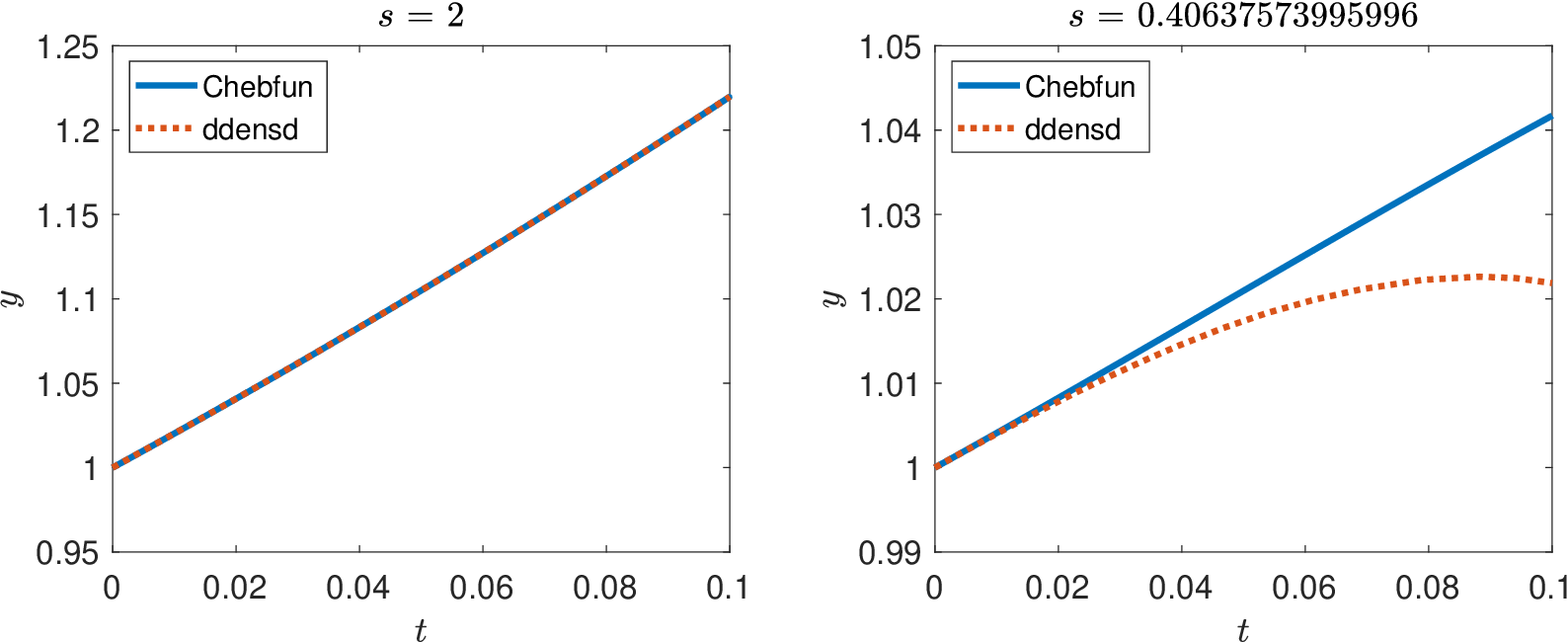}
\end{minipage}
\caption{Top: Chebfun code to solve the neutral-type DDE~(\ref{eqn:chebxample2}). The Chebfun function \lstinline{deriv(y,t/2)} evaluates $y'$ at $t/2$. Bottom left: Agreement of Chebfun and \lstinline{ddesnd} solutions when $s = 2$. Bottom right: Disagreement of Chebfun and \lstinline{ddesnd} solutions when $s =W(-2e^{-2})$.}\label{fig:chebxample2}%
\end{figure}

Our next example is a nonlinear DDE of neutral type:
\begin{equation}\label{eqn:chebxample2}
y'(t) = 2\cos(2t)y(t/2)^{2\cos(t)} + \log(y'(t/2))
              - \log(2\cos(t)) - \sin(t), \quad y(0) = 1.
\end{equation}
This DDE was introduced by~\cite{Jackiewicz1984} and appears in the MATLAB example \lstinline{ddex5} where it is solved using the built-in MATLAB routine for neutral DDEs, \lstinline{ddensd}~\citep{ddensd}.  By substituting $y(0) = 1$ and $y'(0) = s$ to~(\ref{eqn:chebxample2}), one can show that the DDE has two consistent initial slopes (i.e., choices of $y'(0)$ so that~(\ref{eqn:chebxample2}) is satisfied at $t = 0$), namely $s = 2$ and $s = -W(-2e^{-2}) = 0.406\ldots$, where here $W$ is the Lambert W function, each corresponding to a different solution~\citep{Shampine2008}. Chebfun code for computing both solutions is given at the top of Figure~\ref{fig:chebxample2}. When $s = 2$, the DDE has an explicit solution $y = e^{\sin(2t)}$, and both the Chebfun and \lstinline{ddensd} show good agreement with this (with relative errors of around $10^{-14}$ and $10^{-5}$, respectively), as shown in the left panel of Figure~\ref{fig:chebxample2}. For $s = -W(-2e^{-2})$ there is no known closed form solution. In the right panel of Figure~\ref{fig:chebxample2} we see that   the numerical solutions computed by Chebfun and \lstinline{ddesnd} differ considerably, even over a relatively short time interval. This author believes the cause of this is that the DDE is highly ill-conditioned in this case, but that this is mitigated by the fully implicit (i.e., global) nature and high accuracy of the interpolation approach implemented in Chebfun. We defer further investigation to future work.
\subsection{\textbf{\em Chebfun example 3.} (functional, state-dependent)}
Here we simply show how Example~\ref{example:9} can be implemented as a three-line code in Chebfun. The solution is depicted in Figure~\ref{fig:example6}.

	\hspace*{30pt}\begin{minipage}[]{.9\textwidth}
		\begin{lstlisting}
N = chebop(@(t,y) diff(y) + y(y), [0 1]);
N.lbc = 1;
y = N\0;
		\end{lstlisting}
	\end{minipage}%
\subsection{\textbf{\em Chebfun example 4.} (unknown delay)}
\begin{figure}[!t]
\hspace*{13pt}\begin{minipage}[]{.69\textwidth}
\begin{lstlisting}
N = chebop(@(t,y,p) diff(y) + y + y(p*t) - exp(-t/2), [0 1]);
N.lbc = @(y,p) y - 1; 
N.rbc = @(y,p) y - 0.25;
yp = N\0;
plot(yp{1}), title(['p = ' num2str(yp{2}, 10)])
\end{lstlisting}
\end{minipage}\hspace*{3.75pt}
\begin{minipage}[]{.2\textwidth}
\includegraphics[width=103pt]{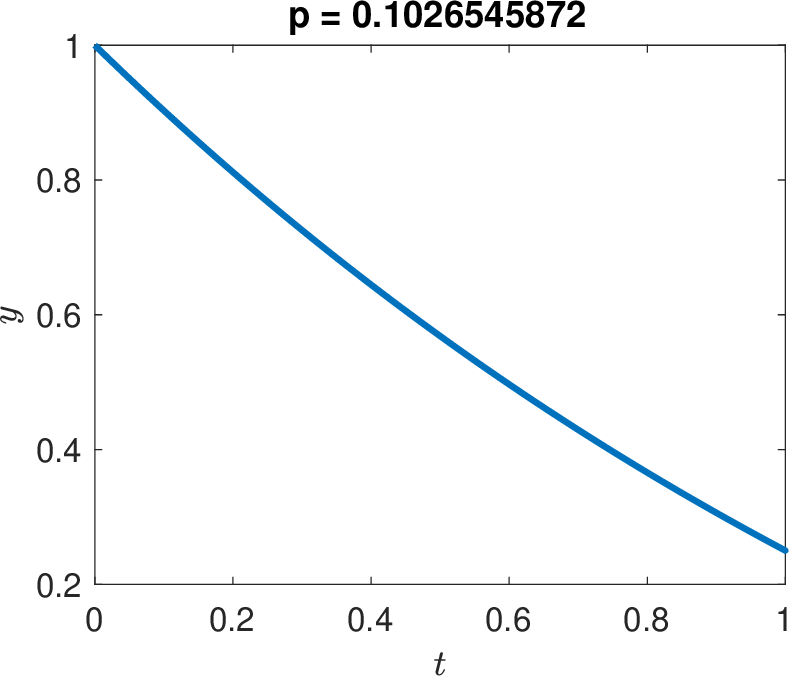}
\end{minipage}
\caption{Left: Chebfun code to solve the unknown delay problem~(\ref{eqn:chebxample4}). Right: Computed solution.}\label{fig:chebxample4} 
\end{figure}

Here we revisit Example~\ref{example:2} from Section~\ref{subsec:proportional}. However, rather than specify the proportional delay, $p$, we consider this an unknown which must be solved in order to satisfy the two-point DDE boundary-value problem  
\begin{equation}\label{eqn:chebxample4}
y'(t) = - y(t) - y(pt) + e^{-t/2}, \qquad y(0) = 1, \ y(1) = 0.25.
\end{equation}
Here we employ Chebop's functionality for solving systems of equations containing unknown parameters without the need to introduce extra equations into the system. With $p$ unknown, the problem is nonlinear, and Chebfun's Newton--Raphson iteration with default initial guess and parameters converges in around 5 iterations. The computed solution (represented by a degree 10 chebfun) and the determined value of $p$ are shown in Figure~\ref{fig:chebxample4}.

\subsection{\textbf{\em Chebfun example 5.} (eigenvalue problem)}
As our final example in this section we consider a delay differential eigenvalue problem (DDE-EVP).\footnote{The most common eigenvalue problem associated with functional or delay differential equations is that of obtaining the characteristic roots of the Jacobian at critical points to determine their stability properties. This is not what we are doing here.}
Chebfun solves differential eigenvalue problems by computing an increasingly large Chebyshev spectral discretisation of the differential operator and looking for the `smoothest' eigenvalues, as described in~\cite[Section 5]{Driscoll2016}.  With the barycentric resampling approach, this is easily adapted to the delay or functional case. Here we consider the second-order DDE-EVP
\begin{equation}\label{eqn:chebxample5}
y''(t) = \lambda y(t/2), \qquad y(0) = y(1) = 0. 
\end{equation}
Chebfun code is given in Figure~\ref{fig:chebxample4}, in addition to the computed eigenvalues and eigenfunctions.
\begin{figure}[!t]
\hspace*{13pt}\begin{minipage}[]{.47\textwidth}
\small
\begin{lstlisting}
A = chebop(@(y) diff(y,2), [0 1]);
A.bc = 'dirichlet';
B = chebop(@(t,y) -y(t/2), [0 1]);
[V, D] = eigs(A, B);
fprintf('%10.10g\n', diag(D))
plot(V)
\end{lstlisting}\vspace*{.5em}
\end{minipage}
\hspace*{5pt}\begin{minipage}[]{.13\textwidth}
\small
\begin{lstlisting}
13.05485001
169.7286494
1398.543635
9480.135715
57516.69512
326708.2900
\end{lstlisting}\vspace*{.5em}
\end{minipage}\hspace*{12.75pt}
\begin{minipage}[]{.3\textwidth}\hfill
\includegraphics[width=103pt]{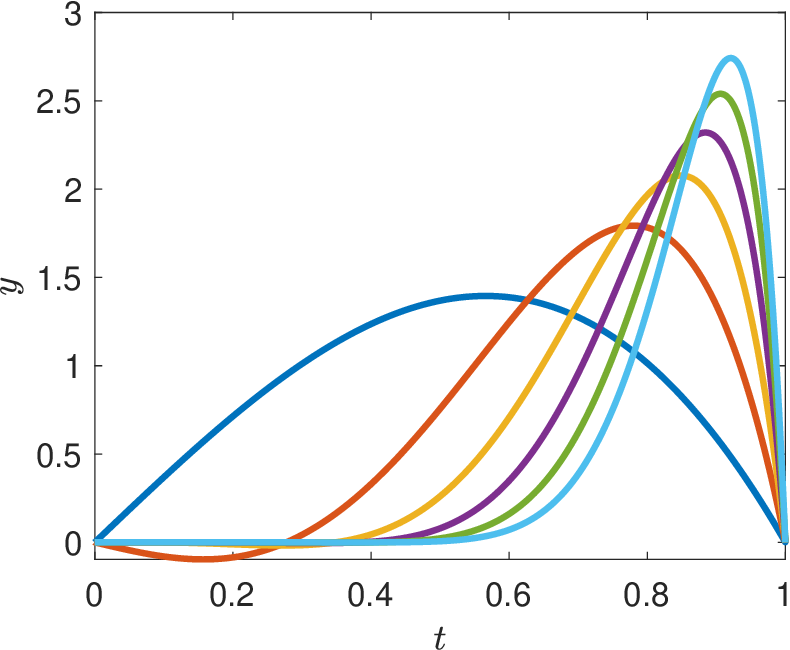}   
\end{minipage}
\caption{Left: Chebfun code to solve to solve the DDE-EVP~(\ref{eqn:chebxample5}). Centre: Computed eigenvalues: Right: Computed eigenfunctions. }\label{fig:chebxample5} 
\end{figure}  


\section{Extensions}\label{sec:ext}%
\subsection{Laguerre and Hermite spectral methods}
Laguerre and Hermite spectral collocation methods, which employ weighted-polynomial bases, $\psi_L(t) = e^{-bt/2}p(t)$ and $\psi_H(t) = e^{-bt^2/2}p(t)$, $b > 0$, have been widely used in solving ODEs and PDEs on semi-infinite and infinite domains, respectively~\citep[Chapter 17]{boyd2001}. The approach described in Sections~\ref{sec:main} and~\ref{sec:functional} above can be naturally extended to solve functional and delay differential equations using these methods as follows. Consider, for example, the Laguerre case, where the weighted-polynomial interpolant may be represented as $\psi_L(t) = e^{-bt/2}{P(t;\bm{t})(e^{b\bm{t}/2}\cdot \bm{\psi})}$, where $\bm{t} = [t_1, \ldots, t_n]^\top$ is an $n$-point Gauss--Laguerre (or Gauss--Laguerre--Lobatto) grid, $P(t;\bm{t})$ is the Lagrange interpolation quasimatrix~(\ref{eqn:quasimatrix}), and $\bm{\psi} = \psi_L(\bm{t})$. It follows that $\psi_L(\tau) = e^{-b\tau/2}{P(\tau;\bm{t})(e^{b\bm{t}/2}\cdot \bm{\psi})}$
and the interpolation operator can be discretised as 
$y(\bm{\tau}) \approx P_L(\bm{\tau};\bm{t})y(\bm{t}),$ where
\begin{equation} P_L(\bm{\tau};\bm{t}) = 
{\rm diag}(e^{-b\bm{\tau}/2})P(\bm{\tau};\bm{t}){\rm diag}(e^{b\bm{t}/2}),\end{equation}
and $P(\bm{\tau};\bm{t})$ is precisely the barycentric resampling matrix from Section~\ref{subsec:bary}. Corresponding differentiation matrices can be constructed using, for example, DMSUITE~\citep{DMSUITE}, and combined with the above as in previous sections. 
This approach has recently been employed by~\cite{HTW2023} for the numerical investigation of the nonlinear functional differential equation
\begin{equation} u^{\prime}(t) +  u(t) =  u^2(\alpha t), \quad \alpha > 0, \  t > 0, \qquad u(0) = 1, \end{equation}
which yields the nonlinear system 
\begin{equation} (D+I)\bm{u} - P_L(\alpha\bm{t};\bm{t})(\bm{u}\cdot\bm{u}) = 0.\end{equation}
\subsection{Fourier spectral methods}

Thus far we have focussed on spectral methods based on polynomial interpolation. However, similar ideas can be applied when solving PDDEs using trigonometric interpolants. In particular, the degree $n$ trigonometric polynomial interpolant of an arbitrary function at $2n$ equidistant nodes, $\theta_k = k\pi/n$, $k = 0,\ldots, 2n-1$ on the interval $[0, 2\pi]$ may be expressed in barycentric form as~\citep{henrici79}
\begin{equation}
t_{n}(\theta) = \sum_{k=0}^{2n-1}(-1)^k\cot\frac{\theta-\theta_k}{2}y_k\bigg/\sum_{k=0}^{2n-1}(-1)^k\cot\frac{\theta-\theta_k}{2}, 
\end{equation}
and a trigonometric barycentric resampling matrix constructed as in the left-hand side of Figure~\ref{fig:barytrig}. This matrix can be combined with the trigonometric differentiation matrix---see~\cite[Chapter 1]{trefethen2000}---to solve DDEs on periodic domains in a similar manner to that described in Section~\ref{sec:main}. An advantage here is that the assumption of periodicity prevents the introduction of discontinuities from non-smooth initial data. 
\begin{example}
Let us first consider a simple linear PDDE, namely 
\begin{equation}\label{eq:periodicexample}
 u''(t)  + \sin(t)u'(t-\pi/\sqrt{2})  + \cos(t)u(t-\pi/2) = 1, \qquad t\in[0, 2\pi],
\end{equation}
with $u(t)$ $2\pi$-periodic. Code to solve this periodic problem is show in the right-hand side of Table~\ref{fig:barytrig}. Figure~\ref{fig:example_trig} shows the computed solution with $n=32$ (left) and the convergence (centre) of the approach as the discretisation size is increased. In this case, we do not have a closed form expression for this solution and so take the $n = 32$ solution as a proxy.  The Fourier coefficients of the $n = 32$ solution are shown in the rightmost panel.
As in the case of the non-periodic problems, we observe geometric convergence in the number of degrees of freedom, in this case down to an accuracy of around $10^{-12}$.

\begin{table}[t!]
\begin{tabularx}{\textwidth}{X|l} 
\end{tabularx}
\hspace*{20pt}\begin{minipage}[]{.43\textwidth}
\begin{lstlisting}[morekeywords={end}]
function P = Barymattrig(z, t)
    T = t(end)-2*t(1)+t(2);
    P = cot((z-t.')*pi/T);
    P(:,2:2:end) = -P(:,2:2:end);
    P = P./sum(P,2);
    P(isnan(P)) = 1;
end
\end{lstlisting}
\end{minipage}
\begin{minipage}[]{.55\textwidth}  
\begin{lstlisting}
t = linspace(0, 2*pi, n+1)'; t(end) = [];
D1 = diffmat(n, 1, 'periodic', [0 2*pi]);
D2 = diffmat(n, 2, 'periodic', [0 2*pi]);
P1 = Barymattrig(t-pi/sqrt(2),t);
P2 = Barymattrig(t-pi/2,t);
A = D2+diag(sin(t))*P1*D1+diag(cos(t))*P2;
u = A\ones(2*n,1);  
\end{lstlisting} 
\end{minipage}
\caption{\footnotesize \rm Left: {MATLAB} code for constructing the trigonometric barycentric resampling matrix for an even number of interpolation nodes \lstinline{t}. For an odd number of points, simply replace \lstinline{cot} with \lstinline{csc}. Right: MATLAB code for solving~\eqref{eq:periodicexample}. Here we use \lstinline{diffmat} from Chebfun to compute the first- and second-order trigonometric differentiation matrices, \lstinline{D1} and \lstinline{D2}.}\label{fig:barytrig}
\end{table}

\begin{figure}[!t] 
\begin{minipage}[]{\textwidth}
	\centering
\includegraphics[height=105pt]{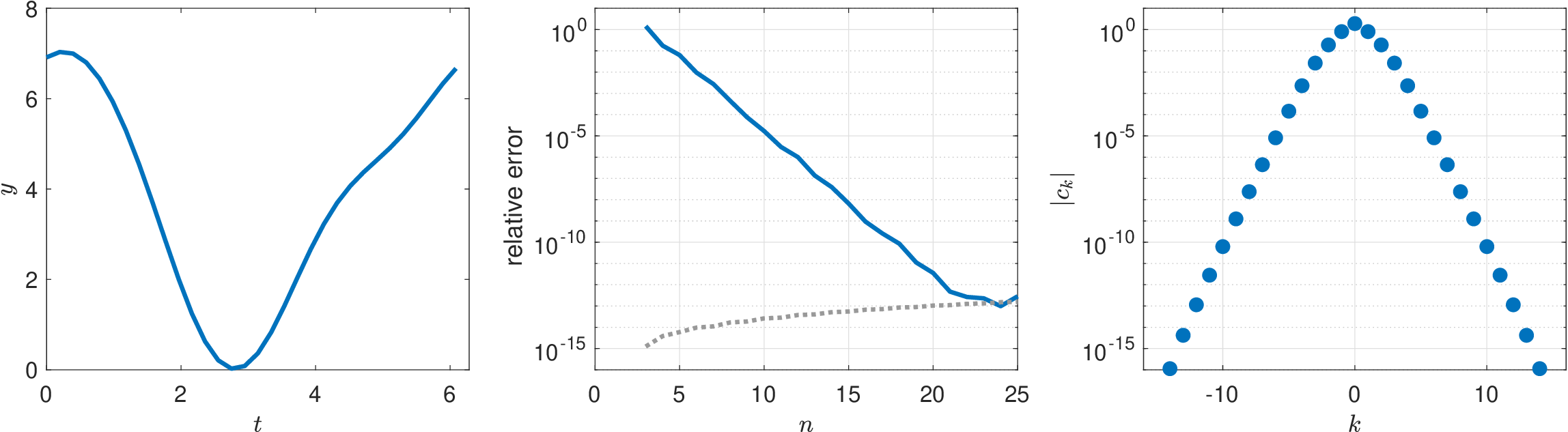}%
\end{minipage} 
\caption{Left: Solution to the periodic DDE~(\ref{eq:periodicexample}) using Fourier collocation on a 32-point uniform grid. Centre: Geometric convergence as $n$ is increased. Right: Fourier coefficients of the computed $n = 32$ solution, exhibiting geometric decay.}\label{fig:example_trig}  
\end{figure} 
\end{example} 

\begin{example}\label{example:11}
A more interesting application is the use of these methods in the accurate computation of limit cycles for autonomous DDEs, which can subsequently be used in the computation of the relevant Floquet multipliers~\citep{Luzyanina2002}. 
For example, consider the following Lotka--Volterra predator-prey model with a delayed Holling type-II response: 
\begin{equation}\label{eqn:lotka}  
 \begin{array}{llll}
  x'(t) & = & x(t) - \displaystyle\frac{x(t)^2}{K} - \frac{x(t)y(t-s)}{1+x(t)},\\
  y'(t) & = & -\gamma y(t)  + \delta \displaystyle\frac{x(t)y(t-s)}{1+x(t)},
 \end{array}, \qquad t > 0.
\end{equation}
For certain parameter choices, one can show that the solution exhibits a stable limit cycle. Here we aim to compute this limit cycle accurately using a Fourier spectral method. The procedure is as follows:
\begin{itemize}
 \item[(a)] For some suitable initial condition $(x_0,y_0)$ in the neighbourhood of the equilibrium point $(x_e,y_e) = (\tfrac{\alpha}{1-\alpha}, (1-\tfrac{x_e}{K})(1+x_e))$, solve~(\ref{eqn:lotka}) for some time using an appropriate time-integration routine (such as \lstinline{dde23} in MATLAB) to obtain an approximation to the limit cycle and its unknown period, $T$. 
 \item[(b)] Rescale~(\ref{eqn:lotka}) to a periodic problem on $[0, 1]$ by an affine transformation of $t$ in the unknown frequency, $1/T$. 
 \item[(c)] Use the Fourier spectral collocation method with barycentric resampling to implement the delay and estimates of the limit cycle and its period from (b) as an initial guess, and solve the periodic nonlinear DDE for $x$, $y$, and $T$.
 \item[(d)] Finally, rescale time to obtain the limit cycle solution on $[0, T]$.
\end{itemize}
A similar approach is employed by~\citep{Borgioli2020} in the context of stability analysis of linear time-periodic delayed dynamical systems with piecewise-smooth coefficients. In Figure~\ref{fig:example_lotka}  we show the result of this procedure when solving~(\ref{eqn:lotka}) with $K = 7/5$, $\gamma = 2/15$, $\delta = 1$, and $s = 1$, for which the period of the limit cycle is found to be $T = 30.83847284$. Checking the residual of the solution with Chebfun gives a forward error of around $6\times10^{-14}$. The solution is depsicted in Figure~\ref{fig:example_lotka}.
\begin{figure}[!t] 
\begin{minipage}[]{\textwidth}
\centering	\includegraphics[height=105pt]{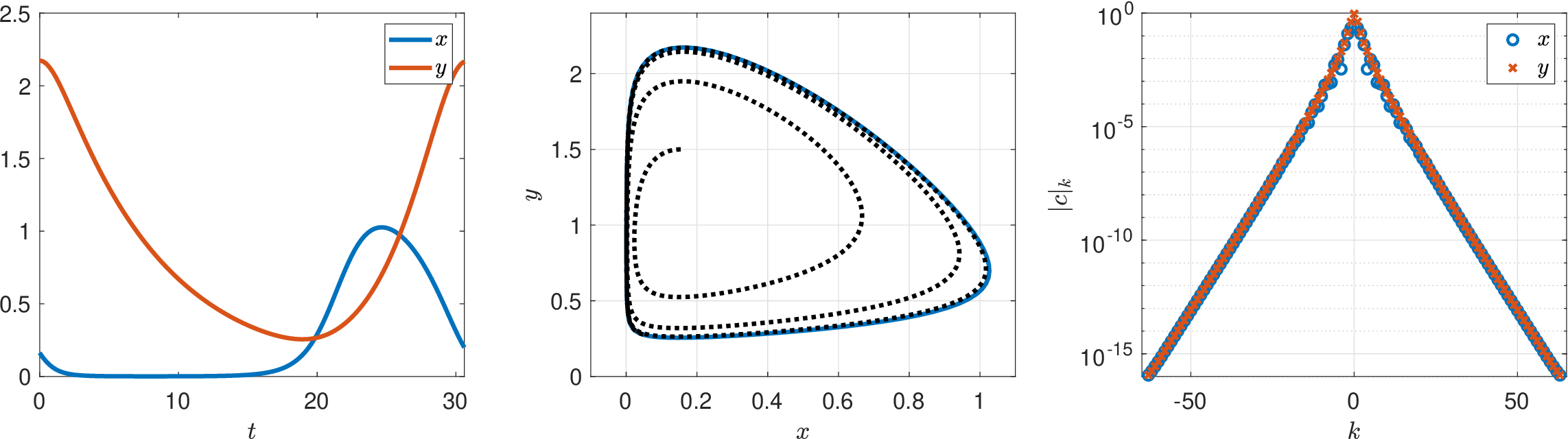}
\end{minipage} 
\caption{Left: Limit cycle solution of the delayed Lotka--Volterra DDE~(\ref{eqn:lotka}), computed using a  129-point grid. Centre: Initial trajectory computed by \lstinline{dde23} (dashed) and limit cycle computed by the Fourier pseudospectral method (solid). Right: Fourier coefficients of the computed  limit cycle solution, demonstrating geometric decay. }\label{fig:example_lotka}  
\end{figure} 
\end{example}

\subsection{Rational spectral methods}

In addition to polynomial and trigonometric interpolants, the approach described in Sections~\ref{sec:main} and~\ref{sec:functional} can be applied to spectral methods based on rational interpolants defined by their barycentric weights, such as those described in~\citep{berrut2014} and~\citep{KTE} for uniform or near-uniformly-spaced nodes, respectively, or those in~\citep{HaleTee} for nodes clustered to resolve local features in the solution. Furthermore, the advent of the AAA algorithm for rational approximation has sparked a renewed interest in numerical methods based on rational functions~\citep{AAA,AAA5}. The AAA algorithm involves a barycentric representation, and as such any spectral collocation method based on AAA ideas for BVP ODEs could be naturally extended to DDEs and FDEs.     

\subsection{Functional equations}

The focus of this work has been on functional {\em differential} equations, but the area of functional {\em equations} (FEs) is in itself of interest, and there is scope for applying the  methods described in this manuscript in such a setting.

\begin{example}  
Consider the Schröder-type functional equation~\citep{Schroeder}
\begin{equation}\label{eq:functional_eqn}
u(f(t)) = \lambda u(t), \qquad f(t) = \lambda\sin(t).  
\end{equation}
For $0 < \lambda < 1$, the theorem of~\cite{Koenigs}  states that there is a unique solution holomorphic  in a neighbourhood of the origin satisfying $u(0) = 0, u'(0) = 1$. Here we approximate this solution on a finite interval by interpolating at scaled  Chebyshev points and employing the barycentric resampling matrix as described in previous sections, yielding a system $\big(P(f(\bm{t});\bm{t})-\lambda I\big)\bm{u} = 0$. This can be augmented with the boundary constraint, $D_{1,:}\bm{u} = 1$, and then solved for $\bm{u} \approx u(\bm{t})$. 
Figure~\ref{fig:functional_eqn} shows code for implementing the above on the domain $[0, \pi]$ with $\lambda = \tfrac12$ and compares the result with the iterated composition solution $u(t) = \lim_{n\rightarrow\infty}\lambda ^{-n}f^{\circ n}(t)$, where $f^{\circ n}(t) := \underbrace{f(f(f(\ldots)))}_{n \text{ times}}.$ As with the FDE examples of previous sections, we observe geometric convergence as the number of interpolation points is increased. The interpolation approach has the advantage that it may be applied to more general FEs where iterated solutions are not readily available.

\begin{figure}[!t]
\hspace*{15pt}\begin{minipage}[]{.42\textwidth}
\begin{lstlisting}
lam = 0.5; f = @(t) lam*sin(t);  
n = 30; [t,~,v] = chebpts(n,[0,pi]);
P = Barymat(f(t), t, v); I = eye(n)
A = [D(1,:) ; 
     P(2:end,:) - lam*I(2:end,:);                          
rhs = [1 ; zeros(n-1,1)];
u = A\rhs;                      
\end{lstlisting}
\end{minipage} 
\begin{minipage}[]{.55\textwidth}\vspace*{8pt}
\includegraphics[height=85pt]{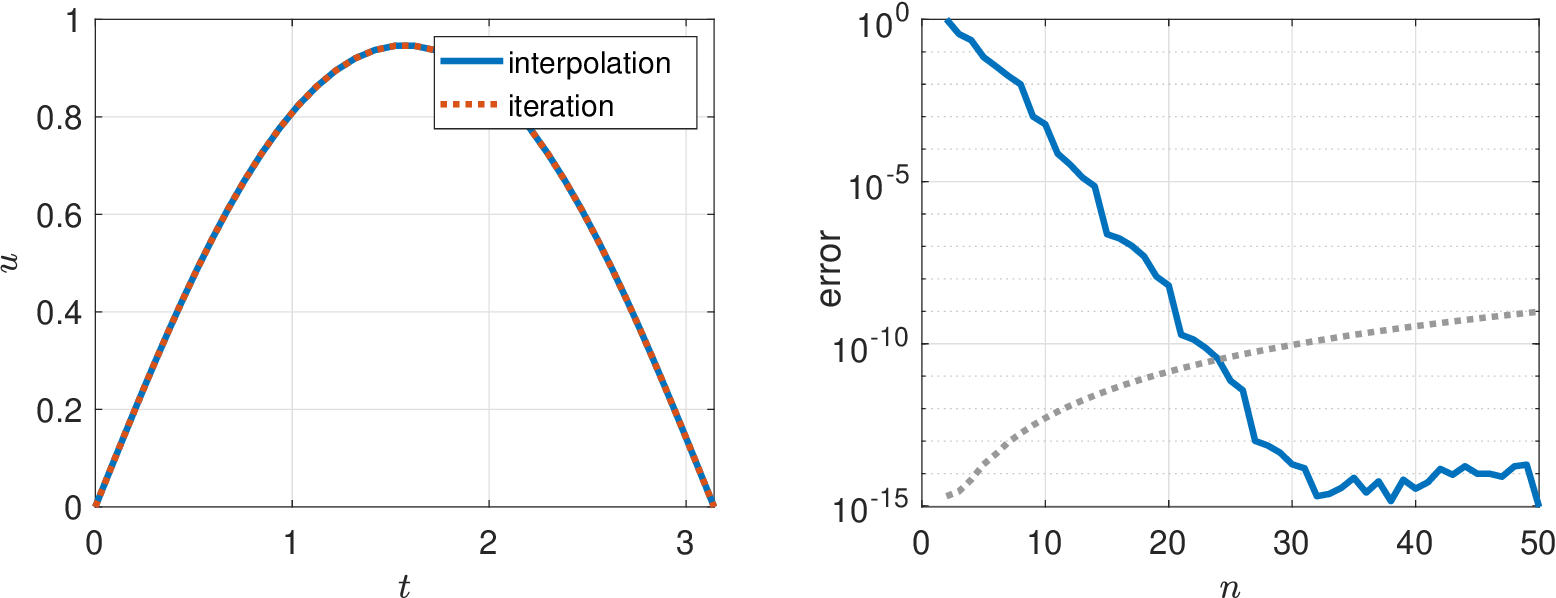}  
\end{minipage}
\caption{Left: MATLAB code to solve the functional equation~(\ref{eq:functional_eqn}) using the barycentric resampling approach.  Centre: Agreement of the computed solutions using interpolation and iterated composition. Right: Geometric convergence of the interpolation solution to the iterated solution as the discretisation size is increased.}\label{fig:functional_eqn}
\end{figure}

\end{example}


\section{Conclusion}\label{sec:conc}%
The barycentric resampling matrix was used to create a framework enabling the solution of various delay and functional differential equations---including pantograph, discrete, continuous, and state-dependent delay---via Chebyshev spectral collocation. Furthermore, this framework was shown to be applicable to neutral-type equations, higher-order DDEs/FDEs, and DDE-EVPs. The close connection to standard Chebyshev spectral methods for ODEs enables the ready implementation of the method to existing code bases, as was demonstrated here in the case of Chebfun/Chebop. Numerous examples were presented, demonstrating both rapid (typically geometric) convergence as well as convenient implementation. Furthermore, the method was shown to be readily extendable to methods based on other classical orthogonal polynomials, trigonometric polynomials, or rational functions. We close with a few additional observations and comments.

For IVP-DDES, it is unlikely that the proposed method will competitive with standard time-stepping--plus--interpolation techniques, for the same reason that spectral methods are rarely used for IVP-ODES (i.e., waveform relaxation methods). However, the example in Section~\ref{sec:chebfunexample2} demonstrates that there may still be some stability benefits to such a global approach for some DDE-IVP problems. For DDE-BVPS, PDDEs, DDE-EVPs, or FDEs that cannot be solved with a method of steps-type approach, there are likely benefits to a spectral collocation method in terms of efficiency and/or accuracy when compared to other approaches---be they shooting-type methods or lower-order collocation approaches---particularly when a high degree of accuracy is required.

In the examples above, we dedicated little discussion to the computational considerations of solving the linear systems of equations resulting from the spectral collocation. For the most part, this is because our focus was
, in spectral element terminology, on the "high $p$, large $h$" regime (i.e., high degree polynomials on a coarse mesh), giving rise to dense matrices and $\mathcal{O}(n^3)$ solves when using direct methods. However, there is perhaps still structure than can be exploited in some cases, particularly as the number of subintervals is increased. Here we discuss some such possibilities. 


In Example~\ref{example:4} we mentioned that the resulting block lower triangular (BLT) discretisation could be efficiently solved via the method of steps for an IVP-DDE. For a BVP-DDE, the augmented boundary conditions will typically result in an `almost block lower triangular' system (i.e., BLT with small number of dense rows). The BLT structure can then be exploited by an appropriate Schur-complement factorisation or the Woodbury formula. These ideas may also be applicable in the FDE case when $\tau(t)\approx t$, leading to almost block-banded systems that could be efficiently solved in a similar way. Another situation that might result in better efficiency for some autonomous problems with periodic boundary conditions, discretised on uniform intervals with the same discretisation size. In such cases,  the resulting linear system has a block circulant structure, which can be exploited~\citep{Mazancourt1983}. However, a global trigonometric interpolant, as in Example~\ref{example:11}, may be better suited to such cases.

Rather than direct methods for solving the linear systems, we might consider iterative approaches. In particular, by utilising the connection to the discrete cosine transform (DCT), Chebyshev differentiation matrices can be applied in $\mathcal{O}(n\log(n))$ operations, making methods like BiCG or GMRES (with suitable preconditioners) attractive~\citep{canuto1985}. In general, the interpolation operator/barycentric resampling matrix requires $\mathcal{O}(n^2)$ operations to apply directly, making an iterative solver less attractive in the DDE/FDE setting. However, this can perhaps be mitigated by the use of a non-uniform DCT~\citep{Potts2003} or fast multipole methods~\citep{dutt1996}.

In the coded examples given, boundary bordering was used to implement boundary conditions. If rectangular projection is used, a la~\citep{Driscoll2016}, then one could explicitly construct the discrete evaluation operators to return values on the projected grid, thus avoiding an additional matrix multiplication. Formulas for the explicit construction of the rectangular differentiation matrices are also known~\citep{XuHale2016}.  

An alternative to spectral collocation methods are spectral methods based on coefficient expansions (e.g., Galerkin or tau methods). The ultraspherical spectral method is such an example, which gives rise to sparse (almost-banded) matrices when applied to ODEs with polynomial coefficients~\citep{olver2013}. One could introduce an interpolation operator to these discretisations, similarly to that in Section 3. This would destroy the sparseness of the discretisation, but other beneficial properties, such as improved conditioning, would likely be maintained, making such an approach attractive. 

Finally, rather than the ``high $p$, large $h$'' approach, one could instead interpolate with lower degree interpolants on a finer mesh, trading exponential convergence with dense matrices for algebraic convergence with sparse matrices. The latter is essentially the approach taken by DDE-BIFTOOL~\citep{BIFTOOLMANUAL}. A detailed comparison of the two approaches is beyond the scope of this paper, but the most efficient approach will likely depend on both the accuracy requirements and on the problem at hand. For smooth problems, such as Example~\ref{example:9} or PDDEs with smooth coefficients, it is typically difficult to beat the geometric convergence of spectral methods by introducing artificial discontinuities in the discretisation~\citep[p.\ {\em x}]{trefethen2000}. For non-smooth problems with multiple discontinuities present in the solution itself, for example, a non-periodic discrete delay DDE-BVP with $p \ll T$, a lower-order discretisation will likely be more efficient. 

\begin{example}\label{example:13}
To demonstrate the former case, we close with one final example, namely computing a limit cycle of the delayed Logistic equation
\begin{equation}\label{eqn:example13}
	y'(t) = (\lambda - y(t - 1))u(t), \qquad y(t+T) = y(t), \quad \forall t.
\end{equation}
This example is considered in~\citep{engelborghs2001}, where a piecewise 3-point Gauss--Legendre discretisation was used to calculate a stable limit cycle with period $T = 4.0964$ for $\lambda  = 1.7$. Using $200$ uniform subintervals, an accuracy of around $5\times10^{-9}$ throughout the interval of integration is achieved, equating to 600 degrees of freedom. As shown in Figure~\ref{fig:final} below, the same accuracy can be achieved using around 50 degrees of freedom with a global polynomial interpolant and only 25 degrees of freedom with a global trigonometric polynomial interpolant.\footnote{In~\citep{engelborghs2001} the error at the interpolation nodes is shown to be lower, about $10^{-11}$, due to super convergence. However, the number of degrees of freedom required to achieve this accuracy in the spectral methods is only increased marginally, due to the geometric convergence.} Of course, degrees of freedom are not the whole story, as one may more readily take advantage of sparsity and or adaptivity in the piecewise discretisation. However, this example suggests that there may be some benefit to global discretisations for such problems. We hope to pursue further investigations in future work.

\begin{figure}[!t]
	\begin{minipage}[]{\textwidth}
		\centering
		\includegraphics[height=105pt]{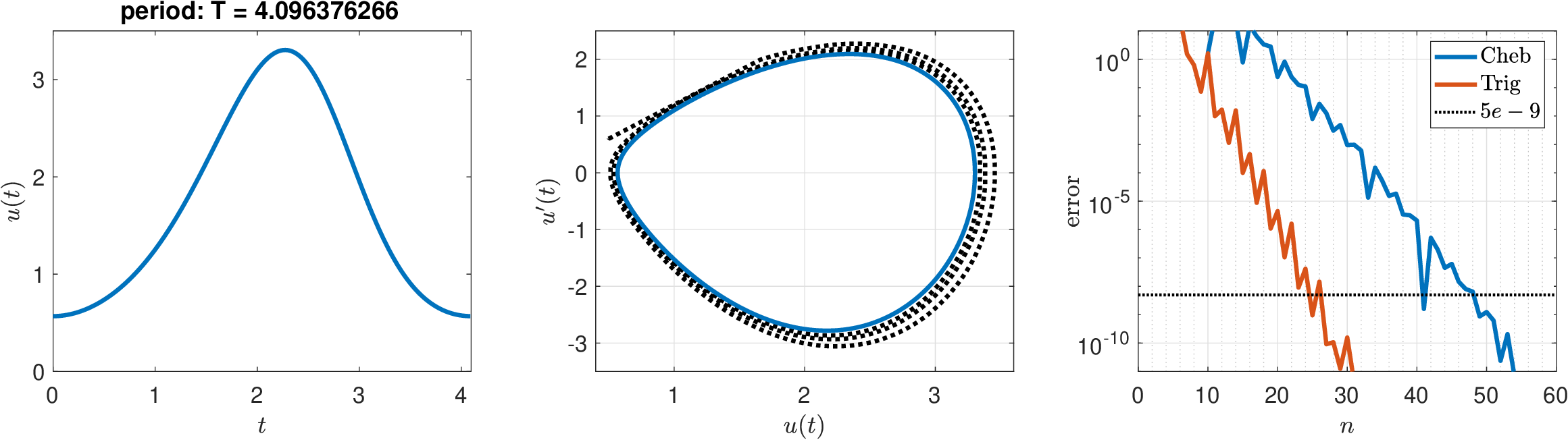}  
	\end{minipage}
	\caption{Left: Periodic solution computed using the spectral collocation method. An initial approximation to the periodic solution is obtained using \lstinline{dde23}, with the initial guess $T = 4$ and assuming $y(t\le0) = 0.5$. This is then refined using the approach described in Example~\ref{example:11}. Centre: Computed limit cycle in the $(y,y')$ plane (solid) and the initial trajectory from \lstinline{dde23} (dashed). Right: Geometric convergence of the spectral collocation approach using both Chebyshev and trigonometric interpolants. (The $n=70$ and $n=40$ solutions are respectively used as a proxy for the unknown analytic solution.)}\label{fig:final}
\end{figure} 

\end{example}

%

\section*{Acknowledgements}
The author is grateful to Andr\'e Weideman and Toby Driscoll for useful discussions and to the anonymous referees for their detailed and constructive feedback.
%
%

\bibliographystyle{abbrvnat}
\bibliography{mybib}

\end{document}